%% file: main.tex
\documentclass[
               DIV=15,
               ]{scrartcl}  

\input{template/packages}
\input{template/pgfplotsettings}

\input{template/commonCommands}
\input{template/titleAndAuthors}
\input{template/journalFooter}

\drawboundingboxtrue    
\drawboundingboxfalse     

\begin{document}     

\normalem          
\maketitle             
    
\vspace{-1.5cm} 
\hrule 
\section*{Abstract}
\input{template/abstract}

\vspace{0.25cm}
\noindent \textit{Keywords:} \input{template/keywords} 
\vspace{0.25cm}
\hrule  

 
\input{introduction/introduction.tex}

\input{problemformulation/problemformulation.tex}

\input{preconditioning/preconditioning.tex}

\input{algorithm/algorithm.tex}

\input{numericalresults/numericalresults.tex}

\input{conclusion/conclusion}

\section*{Acknowledgements} 
\input{template/acknowledgements}
\clearpage
\bibliographystyle{ieeetr}
\bibliography{references.bib}

\end{document}

%% file: template/packages.tex
\usepackage{color} 
\usepackage{psfrag} 
\usepackage{subfig}
\usepackage{hyperref}
\usepackage{cite}
\usepackage{amsmath}
\usepackage{amssymb}
\usepackage{stmaryrd}

\usepackage{xfrac}
\usepackage{float}
\usepackage{graphics}
\usepackage{grffile}
 \usepackage{pstool} 
\usepackage[textsize=tiny]{todonotes}
\usepackage{placeins}
\usepackage[ruled]{algorithm}                   
\usepackage{algorithmicx}                       
\usepackage[noend]{algpseudocode}                              

\usepackage{multirow}
\usepackage{stackengine}

\usepackage[square,comma,numbers,sort&compress]{natbib}
\usepackage[export]{adjustbox}

\usepackage{pgfplots}

\def\fixcheck{%
    \def\tikzifinpicture##1##2{##2}%
}

\usepackage{pgfplotstable}
\usepackage{filecontents}

\usepackage{tikz}
\usepackage{tikzscale}
\usepackage{tikz-3dplot}

\usetikzlibrary{3d}
\usetikzlibrary{spy}
\usetikzlibrary{arrows,shapes}    
\usetikzlibrary{spy}
\usetikzlibrary{backgrounds}  
\usetikzlibrary{decorations}
\usetikzlibrary{decorations.markings}
\usetikzlibrary{positioning}
\usetikzlibrary{patterns}
\usetikzlibrary{calc} 
\usetikzlibrary{fadings}
\usetikzlibrary{trees}
\usetikzlibrary{hobby}
\usetikzlibrary{matrix}
\usepgfplotslibrary{external} 

\usepackage{booktabs}   
\usepackage{makecell} 
\usepackage{colortbl}   

\usepackage{xspace}
\usepackage{color}     
\usepackage{setspace}   
 
\usepackage{soul}
\usepackage{ulem}
\usepackage{currfile}
\usepackage{import}

\usepackage{authoraftertitle} 
\usepackage{authblk} 

\let\OLDthebibliography\thebibliography
\renewcommand\thebibliography[1]{
  \OLDthebibliography{#1}
  \setlength{\parskip}{0pt}
  \setlength{\itemsep}{0pt plus 0.3ex}
}

%% file: template/pgfplotsettings.tex
\pgfplotsset{compat=1.8}

\newcommand{\findmax}[3]{
    \pgfplotstablesort[sort key={#2},sort cmp={float >}]{\sorted}{#1}%
    \pgfplotstablegetelem{0}{#2}\of{\sorted}%
    \let #3=\pgfplotsretval%
}

\pgfplotscreateplotcyclelist{myCycleListBlackAndWhite}
{%
  {black, mark=o},
  {black, mark=square},
  {black, mark=triangle},
  {black, mark=diamond}, 
  {black, mark=pentagon}, 
  {black, mark=*},
  {black, mark=square*},
  {black, mark=triangle*},
  {black, mark=diamond*}, 
  {black, mark=pentagon*}, 
}

\definecolor{darkgreen}{rgb}{0,0.4,0} 
\definecolor{darkbrown}{rgb}{0.5, 0.396, 0.09}

\pgfplotscreateplotcyclelist{myCycleListColor}
{%
  {blue, mark=o},
  {red, mark=square},
  {darkbrown, mark=triangle},
  {darkgreen, mark=diamond},
  {black, mark=pentagon},
  {blue, mark=*},
  {red, mark=square*},
  {darkbrown, mark=triangle*},
  {darkgreen, mark=diamond*},
  {black, mark=pentagon*},
}       

\pgfplotsset{every axis/.append style= 
              {
                font=\small,
                mark size=2,
                line width = 0.1,
                legend style={font=\small, mark size=3, draw=none, fill=none},
                legend cell align=left,
                cycle list name=myCycleListColor,
              }
            } 
            
\pgfplotstableset
{
  every odd row/.style={before row={\rowcolor[gray]{0.9}}},
  every head row/.style=
  {
    before row=
    {%
      \toprule
    },
    after row=\midrule
  },
  every last row/.style=
  {
    after row=\bottomrule
  }
}

\newif\ifdrawboundingbox

\usetikzlibrary{external} 
\tikzset{external/system call={pdflatex \tikzexternalcheckshellescape
-halt-on-error -interaction=batchmode -jobname "\image" "\texsource"}} 
  
\newcommand{\padnum}[2]{%
  \ifnum#1>1 \ifnum#2<10 0\fi
  \ifnum#1>2 \ifnum#2<100 0\fi
  \ifnum#1>3 \ifnum#2<1000 0\fi
  \ifnum#1>4 \ifnum#2<10000 0\fi
  \ifnum#1>5 \ifnum#2<100000 0\fi
  \ifnum#1>6 \ifnum#2<1000000 0\fi
  \ifnum#1>7 \ifnum#2<10000000 0\fi
  \ifnum#1>8 \ifnum#2<100000000 0\fi
  \ifnum#1>9 \ifnum#2<1000000000 0\fi
  \fi\fi\fi\fi\fi\fi\fi\fi\fi
  \expandafter\expandafter\number#2%
}

\newcounter{tmp}
\newcommand{\nextfigurename}[1]%
{%
	\setcounter{tmp}{\thefigure}%
	\addtocounter{tmp}{1}%
	\tikzsetnextfilename{section\thesubsection_#1}%
}%

\newcommand{\nextsubfigurename}[1]%
{%
	\tikzsetnextfilename{section\thesubsection_#1}%
}%
          
\makeatletter
\renewcommand{\todo}[2][]{\tikzexternaldisable\@todo[#1]{#2}\tikzexternalenable}
\makeatother

\makeatletter
\renewcommand{\missingfigure}[2][]{\tikzexternaldisable\@missingfigure[#1]{#2}\tikzexternalenable}
\makeatother

%% file: template/commonCommands.tex
\newcolumntype{C}[1]{>{\centering\arraybackslash}m{#1}}
\newcolumntype{R}[1]{>{\raggedright\arraybackslash}m{#1}}
\newcolumntype{L}[1]{>{\raggedleft\arraybackslash}m{#1}}

\newcommand{\includetikz}[2][]{%
\StrSubstitute{#2}{/}{_}[\temp]%
\tikzsetnextfilename{\temp}%
\includegraphics[#1]{#2.tikz}%
} 

\newcommand{\delete}[1]{\xspace} 
   
\newcommand{\rootDir}{}
\newcommand{\graphDir}{}
\newcommand{\dataDir}{}
\newcommand{\picsDir}{}

\newtheorem{theorem}{Theorem}
\newtheorem{remark}[theorem]{Remark}

%% file: template/titleAndAuthors.tex
\title{ Robust and parallel scalable iterative solutions \\ for large-scale finite cell analyses}

\author[1]{J. N. Jomo\thanks{Corresponding author: john.jomo@tum.de}}
\author[2]{F. de Prenter}
\author[1]{M. Elhaddad}
\author[1,5]{D. D'Angella}
\author[2]{\authorcr C. V. Verhoosel}
\author[1]{S. Kollmannsberger}
\author[3]{J. S. Kirschke}
\author[4]{\authorcr V. N\"ubel}
\author[2]{E. H. van Brummelen}
\author[1,5]{E. Rank}

\affil[1]{Chair for Computation in Engineering, Technische Universit\"at M\"unchen}
\affil[2]{Department of Mechanical Engineering, Eindhoven University of Technology}
\affil[3]{Abteilung f\"ur Neuroradiologie, Klinikum rechts der Isar, Technische Universit\"at M\"unchen}
\affil[4]{Hilti Entwicklungsgesellschaft GmbH}
\affil[5]{Institute for Advanced Study, Technische Universit\"at M\"unchen}

\newcommand{\journal}{Finite Elements in Analysis and Design}

%% file: template/journalFooter.tex
\usepackage{fancyhdr}
\date{}
\fancypagestyle{plain}{%
  \fancyhf{}%
  \fancyfoot[L]{
  \vspace{-1.5cm} 
  \footnotesize{Preprint submitted to \textit{\journal{}}  \hfill
  }} }

%% file: template/abstract.tex
The finite cell method is a highly flexible discretization technique for numerical analysis on domains with complex geometries. By using a non-boundary conforming computational domain that can be easily meshed, automatized computations on a wide range of geometrical models can be performed. Application of the finite cell method, and other immersed methods, to large real-life and industrial problems is often limited due to the conditioning problems associated with these methods. These conditioning problems have caused researchers to resort to direct solution methods, which significantly limit the maximum size of solvable systems. Iterative solvers are better suited for large-scale computations than their direct counterparts due to their lower memory requirements and suitability for parallel computing. These benefits can, however, only be exploited when systems are properly conditioned. In this contribution we present an Additive-Schwarz type preconditioner that enables efficient and parallel scalable iterative solutions of large-scale multi-level $hp$-refined finite cell analyses.

%% file: template/keywords.tex
immersed methods, finite cell method, iterative solvers, preconditioning, \emph{hp}-refinement, parallel computing, high performance computing

%% file: introduction/introduction.tex
\section{Introduction}
\label{sec:introduction}

Immersed finite element methods, most prominently the finite cell method \cite{Parvizian2007,Duster2008,Schillinger2014,Xu2016,Varduhn2016,Duster2017} and cutFEM \cite{Burman2012,Burman2014}, are useful techniques for performing numerical computations on domains with complex geometries. Generating boundary fitted meshes for these problems may be an error-prone, laborious and computationally expensive operation. The general concept of immersed finite element methods is not to directly mesh the complex geometry on which the problem is posed, but to instead embed it in a larger, embedding domain of simple shape that can be trivially discretized. An approximation space is then defined on this structured mesh and restricted to the complex geometry of the problem domain.

The finite cell method (FCM) combines the embedded mesh concept with high-order finite elements to compute highly accurate results with a flexible discretization scheme. The efficiency of immersed methods can be further improved by using them in conjunction with $hp$-refined discretizations. The multi-level $hp$-finite cell method \cite{Zander2015,Zander2016} has shown to be effective in capturing local solution characteristics, such as steep gradients or stress concentrations induced by (fine) geometrical features, e.g.\ \cite{Zander2016a,Bog2017,Elhaddad2017}. An unavoidable aspect of applying this method to large-scale real-life applications is the necessity of solving large linear systems.
The computational cost of using a direct solver does not scale well with the size of the system, however, and becomes increasingly expensive in terms of memory consumption and execution time. Iterative solvers are the preferred method for solving large systems, because their computational cost scales better with the system size and due to their low memory requirement and ability to be easily parallelized in comparison to direct solvers. Convergence of iterative solvers, however, strongly depends on the conditioning of the system, e.g.\ \cite{Saad2003}. Without tailored preconditioning, FCM systems generally show severe ill-conditioning, which practically prohibits application of iterative solvers \cite{SIPIC}. As a result, finite cell computations have been mainly restricted to the use of direct solvers \cite{Duster2008,schillinger_isogeometric_2012,Rank2012,ruess_weakly_2013,Ruess2014,Schillinger2014}.

Several techniques have been developed in order to improve the conditioning of immersed finite elements and enable iterative solution techniques. Fictitious domain stiffness assumes a soft material in the fictitious domain, and thereby performs a volumetric stabilization through the evaluation of integrals in both the problem and fictitious domain. An extensive analysis of this method can be found in \cite{Dauge2015}. For FCM problems in the context of isogeometric analysis, it is shown in \cite{Verhoosel2015,Elfverson2018} that problematic basis functions with small support inside the domain of computation can simply be removed from the system. Ghost penalty stabilization \cite{Burman2010} is another option which weakly couples the degrees of freedom (DOFs) in cut elements to their neighboring elements. This gives better control of the solution on cut elements and improves the conditioning. Application of ghost penalty is customary in methods referred to as cutFEM \cite{Burman2012,Burman2014}. Instead of weakly coupling DOFs in cut elements to neighboring elements, these can also be constrained strongly as originally introduced in \cite{Hoellig2001,Hoellig2005} and applied in e.g.\ \cite{Sanches2011,Rueberg2012,Rueberg2013,Rueberg2016,Badia2018,Marussig2018}. 

In addition to these stabilization methods multiple preconditioning techniques have been proposed for immersed finite element methods and for similar problems in the extended finite element method (X-FEM) \cite{Belytschko1999,Moes1999}. In X-FEM, these preconditioners generally split the matrix in \emph{(i)} non-problematic DOFs, which can be treated by standard preconditioning techniques and \emph{(ii)} problematic DOFs, that are treated separately by local Cholesky decompositions \cite{Bechet2005}, a tailored FETI-type method \cite{Menk2011}, a Schur complement based algebraic multigrid preconditioner \cite{Hiriyur2012} or a Schwarz-type domain decomposition preconditioner \cite{Berger-Vergiat2012,Waisman2013}. The last two references are conceptually similar to the approach applied in this work, but differ in the choice of subdomains as the problematic DOFs in the finite cell method, i.e.\ all DOFs on the boundary, are generally large connected sets. For immersed finite element methods it is proposed in \cite{Lehrenfeld2017} to diagonally scale the problematic DOFs and treat the non-problematic DOFs by an algebraic multigrid procedure, which is demonstrated to be  effective for linear basis functions and restrictions on the manner in which elements can be cut. In \cite{Badia2017} the scaling in a Balancing Domain Decomposition by Constraints preconditioner (BDDC) is customized for cut basis functions, which is shown to result in an effective method for linear bases. A preconditioner that combines a diagonal scaling with a local orthonormalization of the problematic DOFs is developed in \cite{SIPIC}. In \cite{CbAS} it is shown that this orthonormalization is very similar to Additive-Schwarz preconditioning, e.g.\ \cite{BrennerScott}, of the cut elements. The preconditioner developed in \cite{CbAS} is demonstrated to be effective for immersed finite element methods with higher order basis functions, but has so far only been tested on uniform tensor product meshes without local refinements.

In this contribution, we further develop the preconditioning technique presented in \cite{CbAS} through modifications that are essential to allow robust application of iterative solvers for finite cell discretizations involving local $hp$-refinements. The robustness of the preconditioner is improved by stabilizing the local inverses in case the smallest eigenvalues approach machine precision. Furthermore, the efficiency and scalability of the preconditioner setup is investigated and optimized by \textit{(i)} focusing on problematic areas and \textit{(ii)} the use of shared memory and distributed memory parallelization. The developed technique is demonstrated to be able to efficiently solve large FCM systems involving multi-level $hp$-refinement with good parallel scalability using a preconditioned Conjugate Gradient (PCG) solver \cite{Saad2003}. The presented examples show the potential of using FCM in large-scale real-life and industrial applications.

Section 2 introduces the finite cell method and multi-level $hp$-refinement. Section 3 recapitulates the conditioning of the finite cell method and the principal concepts of Additive-Schwarz preconditioning, which is followed by the modifications required for preconditioning multi-level $hp$-refined discretizations. Section 4 presents the details of the parallel implementation and Section 5 investigates and demonstrates the effectiveness of the preconditioner with a PCG solver in several numerical examples.

%% file: problemformulation/problemformulation.tex
\renewcommand{\rootDir}{problemformulation}
\renewcommand{\graphDir}{\rootDir/graphs}
\renewcommand{\dataDir}{\rootDir/data}
\renewcommand{\picsDir}{\rootDir/pics} 

\section{Discretization techniques}
In this section we introduce the discretization methods that form the basis of the proposed techniques, namely the finite cell method in Section \ref{sec::finitecellmethod} and multi-level $hp$-refinement in Section \ref{sec::hp-discretization}\,.

\subsection{The finite cell method}\label{sec::finitecellmethod} 

\begin{figure}[t!]
\begin{center}
\subfloat[Physical domain $\Omega_{\text{phys}}$]
{
 \includegraphics[width=0.31\textwidth]{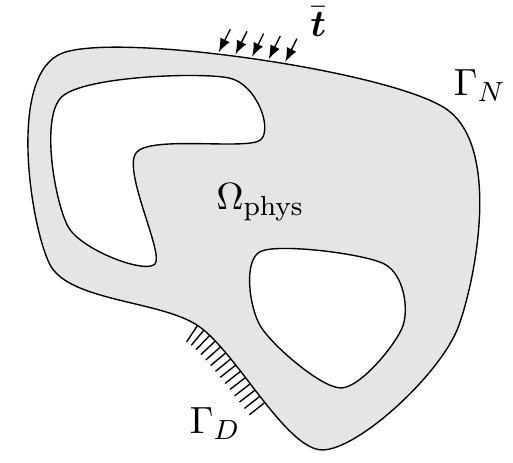}
}
\subfloat[Fictitious extension: $\Omega_{\textrm{phys}} \cup \Omega_{\textrm{fict}}$]
{
 \includegraphics[width=0.31\textwidth]{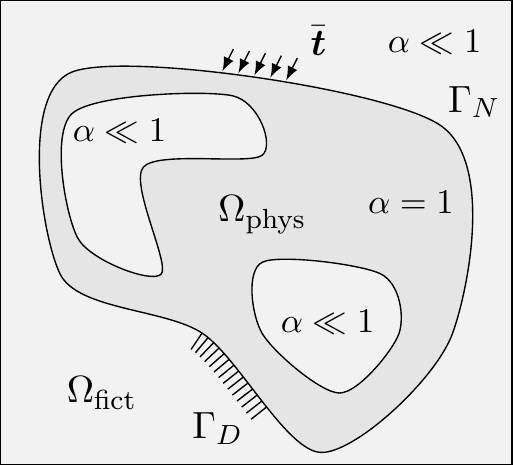}
}
\subfloat[Structured discretization]
{
 \includegraphics[width=0.31\textwidth]{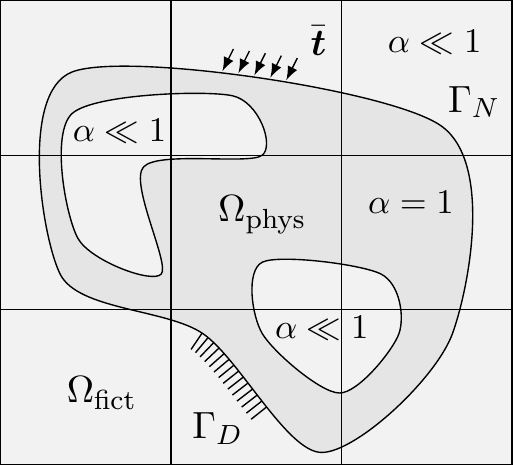}
}
\caption{A schematic representation of the finite cell method \cite{Bog2017}.}
\label{fig::FCMrepresentation}
\end{center}
\end{figure}

The finite cell method (FCM) is an immersed high-order method suitable for numerical analysis on domains of complex shape \cite{Parvizian2007,Duster2008}. As an underlying numerical approach finite elements or spline-based shape functions can be used. In the latter case FCM turns out to be a suitable technique for trimmed isogeometric analysis \cite{Schillinger2011,Kamensky2015}. FCM extends the original problem domain referred to as the physical domain, $\Omega_{\textrm{phys}}$, with a surrounding fictitious domain, $\Omega_{\textrm{fict}}$, forming a computational domain $\Omega$ with a simple structure that can be trivially meshed (Figure~\ref{fig::FCMrepresentation}). Distinction between points in the physical and fictitious domains is made by an indicator function $\alpha$, which assumes the value $\alpha=1$ for points lying in $\Omega_{\textrm{phys}}$ and $0 \leq \alpha \ll 1$ for points outside $\Omega_{\textrm{phys}}$. The parameter $\alpha$ is generally chosen as a small finite constant $0 < \alpha \ll 1$ instead of $\alpha = 0$ to ensure numerical stability, e.g.\ \cite{Dauge2015}.
Since the basis is generally not interpolatory on the physical domain boundary, Dirichlet conditions are weakly imposed using e.g.\ Nitsche's method \cite{Nitsche1971,Embar2010} or the penalty method \cite{Babuska1973}. Because FCM results in conditioning problems of a similar nature with both Nitsche's method and the penalty method, we apply the penalty method in all our examples because of its simplicity and low computational cost. All presented results can be extended to problems with Nitsche boundary conditions without significant changes. 

The weak form in FCM is illustrated for a linear elasticity problem in Eq.\,\eqref{eq::modifiedweakform} where $\boldsymbol{u}$ is the unknown displacement field, $\boldsymbol{v}$ an arbitrary test function, $\boldsymbol{\varepsilon}$ the infinitesimal strain tensor, $\mathbb{C}$ the fourth-order elasticity tensor and $\beta$ the penalty parameter. The vectors $\boldsymbol{b}$, $\boldsymbol{f}$ and $\boldsymbol{g}$ represent the volumetric forces, surface traction and Dirichlet data respectively.  
\begin{equation}\label{eq::modifiedweakform}
	\int \limits_{\Omega} \alpha \, \boldsymbol{\varepsilon}(\boldsymbol{v}) : \mathbb{C} : \boldsymbol{\varepsilon}(\boldsymbol{u}) \, d\Omega \ + \ \int \limits_{\Gamma_D} \beta \, \boldsymbol{v} \cdot \boldsymbol{u} \, d\Gamma   = \int \limits_{\Omega} \alpha \, \boldsymbol{v} \cdot \boldsymbol{b} \, d\Omega \ + \ \int \limits_{\Gamma_{N}} \boldsymbol{v} \cdot \boldsymbol{f} \,  d\Gamma \ + \int \limits_{\Gamma_D} \beta \, \boldsymbol{v} \cdot  \boldsymbol{g} \,  d\Gamma \, . 
\end{equation}
FCM has a large application field due to its ability to produce highly flexible, non-conforming discretizations while maintaining the accuracy and computational efficiency of high-order finite elements. For more details regarding FCM the reader is referred \cite{Parvizian2007,Duster2008} or the review articles \cite{Schillinger2014,Duster2017} that show various applications and recent developments.

\subsection{Multi-level $hp$-discretizations}\label{sec::hp-discretization} 

For problems with non-smooth solutions, a high quality finite element approximation can be obtained with a relatively small number of unknowns by the use of $hp$-refinement schemes. These schemes combine a local increase of the spatial resolution  ($h$-refinement) and an elevation of the element polynomial order ($p$-refinement) in the crucial regions of the computational domain.
Often $hp$-refinement schemes perform spatial refinement by the replacement of large elements having a high approximation error with a set of smaller elements, e.g.~\cite{Solin2004a,Rachowicz2006,Demkowicz2007}. Multi-level $hp$-refinement, introduced in \cite{Zander2015}, performs $hp$-refinement by means of high-order hierarchical overlay meshes. Instead of replacing large coarse elements, multi-level $hp$-refinement superposes coarse elements with a hierarchy of overlay elements in order to better capture solution characteristics in regions of interest. The solution is the sum of the large-scale solution on the coarse base mesh and the fine-scale solutions on the overlay meshes. This change of paradigm from a replacement-driven to a superposition-driven spatial refinement significantly reduces the complexity associated with the implementation of $hp$-refinement schemes. The direct relation of degrees of freedom to topological components (nodes, edges, faces and volumes) in high-order finite elements allows a straightforward enforcement of both the compatibility of the discretization and linear independence of the basis functions in the multi-level $hp$-refinement setting \textit{c.f} \cite{Zander2015,Zander2016} and \cite{Zander2016a,Elhaddad2017} for applications.

\begin{figure}[t]
  \begin{minipage}[c]{0.7\textwidth}
    \begin{minipage}[t]{0.45\textheight}	  
     \begin{center}
       \includegraphics[width=0.95\textwidth]{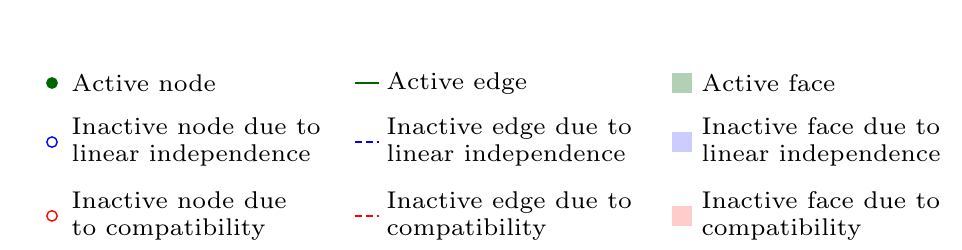}
     \end{center}	
    \end{minipage}	  
    \begin{minipage}[t]{0.5\textheight}
     \centering	    
     \begin{center}
       \subfloat[One-dimensional case]{ \includegraphics[width=0.56\textwidth]{\picsDir/refinement_1d_mlhp_woNumbers} }
       \subfloat[Two-dimensional case]{ \includegraphics[width=0.38\textwidth]{\picsDir/multiLevelhp2d} }%
     \end{center}	
   \end{minipage}	  
  \end{minipage}
  \begin{minipage}[c]{0.28\textwidth}
   \centering	 
   \vspace{4mm} 
   \begin{center}
      \includetikz[width=0.95\textwidth]{\picsDir/multiLevelhp3d-notext}%
    \end{center}
    \footnotesize
    (c) Three-dimensional case
  \end{minipage}
  \caption{Illustration of the multi-level $hp$-refinement scheme with two refinement levels, $k=2$, in different spatial dimensions. The deactivation of specific topological components following a simple rule-set ensures compatibility and linear independence of the basis functions  \cite{Zander2016}.}
  \label{fig::multilevelhp::multilevelhpIdea}%
\end{figure}

Figure \ref{fig::multilevelhp::multilevelhpIdea} illustrates the superposition principle in multi-level $hp$-refinement. The scheme is not restricted to 1-irregular meshes since arbitrary levels of hanging nodes are treated automatically through the use of homogeneous Dirichlet boundary conditions at the boundary of the overlay meshes. This deactivates specific topological components and enforces $C^0$-continuity over the complete refinement hierarchy. $C^{\infty}$-continuity is maintained in unrefined base elements and overlay elements on the highest refinement level. Such elements are collectively referred to as \emph{active or leaf elements}. Numerical integration is performed on these leaf elements.

%% file: preconditioning/preconditioning.tex
\renewcommand{\rootDir}{preconditioning}
\renewcommand{\graphDir}{\rootDir/graphs}
\renewcommand{\dataDir}{\rootDir/data}
\renewcommand{\picsDir}{\rootDir/pics} 

\section{Preconditioning}\label{sec::preconGeneral}

As mentioned in Section~\ref{sec:introduction}\,, iterative solution methods are preferred over direct solvers for large systems due to their lower memory requirements and suitability for parallel computing. It is well known that the convergence speed of iterative methods generally depends on the conditioning of the system, e.g.\ \cite{Saad2003}. The conditioning of FCM systems is generally rather poor, as remarked in e.g. \cite{schillinger_isogeometric_2012,Rank2012,ruess_weakly_2013,Ruess2014,Schillinger2014} and shown in detail in \cite{SIPIC}, which impedes solving large FCM systems. In \cite{CbAS} an Additive-Schwarz scheme is proposed that is efficient for FCM systems on uniform grids, but that does not account for the local refinements in multi-level $hp$-discretizations. Section~\ref{sec::additiveschwarzfcm} revisits the conditioning problems associated with the finite cell method and discusses how Additive-Schwarz preconditioners can alleviate these problems. Section~\ref{sec::preconditioningmlhprefinement} presents the required modifications to the Additive-Schwarz scheme that allow efficient preconditioning of immersed multi-level $hp$-discretizations. 

\subsection{FCM conditioning and Additive-Schwarz preconditioning}\label{sec::additiveschwarzfcm}

This manuscript focuses on problems in linear elasticity, for which FCM yields a symmetric positive definite (SPD) system of the form
\begin{equation}
\boldsymbol{A} \boldsymbol{x} = \boldsymbol{b} \, ,
\end{equation}
with system matrix $\boldsymbol{A}$, solution vector $\boldsymbol{x}$ and right hand side vector $\boldsymbol{b}$. The condition number $\kappa\left(\boldsymbol{A}\right)$ of SPD systems is defined as
\begin{equation}
\kappa \left( \boldsymbol{A} \right) = \large\| \boldsymbol{A} \large\| \, \large\| \boldsymbol{A}^{-1} \large\| = \frac{\lambda_{\textsc{max}}}{\lambda_{\textsc{min}}} = \frac{\max_{\|\boldsymbol{v} \|=1} \boldsymbol{v}^T \boldsymbol{A} \boldsymbol{v}}{\min_{\|\boldsymbol{u} \|=1} \boldsymbol{u}^T \boldsymbol{A} \boldsymbol{u}} \, ,
\end{equation}
with $\lambda_{\textsc{max}}$ and $\lambda_{\textsc{min}}$ denoting the largest and smallest eigenvalue and $\boldsymbol{v}$ and $\boldsymbol{u}$ (eigen)vectors maximizing and minimizing the Rayleigh quotient with $\boldsymbol{A}$. From \cite{SIPIC} we know that FCM systems are generally ill-conditioned because of the occurrence of elements with a very small intersection with the physical domain. Basis functions on such elements can be arbitrarily small in magnitude and/or basis functions that are sufficiently linearly independent on the full element can become almost linearly dependent when the element is cut. The latter is caused firstly by higher order contributions of basis functions, that can become arbitrarily small compared to linear (or lower order) contributions, and secondly by dependencies on a (e.g. horizontal or vertical) parametric dimension, that can become arbitrarily small compared to other parametric dimensions. This reduces a higher order or multivariate function to an almost linear or univariate function. When this occurs for multiple functions, these functions can become almost linearly dependent. The magnitude of a function in the approximation space that is only supported on a small cut element --- which strongly depends on the volume of the intersection between the element and the physical domain --- can therefore be arbitrarily small compared to the magnitude of the vector that represents the function in the isomorphic vector space --- which is independent of this volume. Because of that, FCM systems generally have $\min_{\|\boldsymbol{u} \|=1} \boldsymbol{u}^T \boldsymbol{A} \boldsymbol{u} \ll 1$, such that $\kappa\left(\boldsymbol{A}\right) \gg 1$. When the smallest cut element is assumed to be shape regular and for $\alpha = 0$, it is shown in \cite{SIPIC} that the condition number of second order FCM systems can be estimated by
\begin{equation}
\kappa \left( \boldsymbol{A} \right) \propto C \eta^{-(2p+1-2/d)} \, ,
\label{eq::kappacut}
\end{equation}
with $p$ the polynomial order of the discretization, $d$ the number of dimensions and $\eta$ the relative volume of the intersection between the smallest cut element and the physical domain
\begin{equation}
 \eta = \min_{T \in \mathcal{T}} \frac{|T \cap \Omega_{\text{phys}}|}{|T|} \, , 
\label{eq::volumefraction}
\end{equation}
with $\mathcal{T}$ the set of elements, $|T \cap \Omega_{\text{phys}}|$ the volume of the intersection between element $T$ and the physical domain and $|T|$ the full volume of element $T$. Since the volume of the intersection can become arbitrarily small, the condition number in FCM can be arbitrarily large.

\begin{remark}
It is stipulated that the ill-conditioning of FCM is not only due to basis functions that become very small. This could easily be resolved by Jacobi or Gauss-Seidel preconditioning. Notwithstanding that Jacobi and Gauss-Seidel preconditioners improve the conditioning, these do not repair conditioning problems related to almost linearly dependent basis functions on small cut elements. Therefore these simple preconditioners do not result in robust preconditioning for the finite cell method, as is shown in \cite{SIPIC} and demonstrated further in Section \ref{sec::vertebra}\,.
\end{remark}

\subsubsection*{Additive-Schwarz preconditioning}

Ill-conditioning in FCM can be effectively resolved by the use of Additive-Schwarz preconditioners \cite{CbAS}. A rich literature exists on these preconditioners in a general setting, e.g.\ \cite{Smith1996,Toselli2005}, and specifically for finite element methods, e.g.\ \cite{Ferencz1998,BrennerScott}. The general idea of Additive-Schwarz preconditioners is the construction of a sparse approximation of $\boldsymbol{A}^{-1}$ through the inversion and summation of sub-matrices of matrix $\boldsymbol{A}$ as per the formula
\begin{equation}\label{eq:AdditiveSchwarz}
	\boldsymbol{S} = \sum_{B \in \mathcal{B}} \boldsymbol{R}_B^T\underbrace{\left(\boldsymbol{R}_B\boldsymbol{A}\boldsymbol{R}_B^T\right)^{-1}}_{\boldsymbol{A}_B^{-1}}\boldsymbol{R}_B \, .
\end{equation}
In this formulation $B$ denotes a \emph{block}, which is a set of indices referring to basis functions. The set of all blocks is denoted by $\mathcal{B}$. A basis function $\phi_k$ can be in multiple (i.e.\ overlapping) blocks, but each basis function must be in at least one block for $\boldsymbol{S}$ to be nonsingular. Premultiplying matrix $\boldsymbol{A}$ with the (non-square) restriction matrix $\boldsymbol{R}_B$ restricts $\boldsymbol{A}$ to the row indices in $B$, and postmultiplying matrix $\boldsymbol{A}$ with the restriction matrix $\boldsymbol{R}_B$ restricts $\boldsymbol{A}$ to the column indices in $B$. The inverse of sub-matrix $\boldsymbol{A}_B$ is placed at the indices of block $B$ by pre- and postmultiplying with $\boldsymbol{R}_B^T$ and $\boldsymbol{R}_B$ respectively.

Properly choosing the blocks $\mathcal{B}$ is essential for the effectiveness of Additive-Schwarz preconditioning, because different blocks yield varying preconditioners with different properties. For example one can assign $\mathcal{B}$ such that every function is in a single, separate block. This results in the standard Jacobi preconditioner which is not robust to cut elements, as demonstrated in Section~\ref{sec::vertebra}\,. On the other extreme one can theoretically assign only one block containing all functions such that $\boldsymbol{S} = \boldsymbol{A}^{-1}$. Such a preconditioner is optimal in terms of spectral properties, but is prohibitively expensive as it involves inverting the full system. For FCM, it is useful to assign blocks based on the Additive-Schwarz lemma \cite{Matsokin1985,Lions1988,Smith1996,Toselli2005}:
\begin{equation}\label{eq:AdditiveSchwarzLemma}
\boldsymbol{v}^T \boldsymbol{S}^{-1} \boldsymbol{v} = \min_{\boldsymbol{v}=\sum\limits_{B \in \mathcal{B}} \boldsymbol{R}_B^T \boldsymbol{v}_B} \sum_{B \in \mathcal{B}} \boldsymbol{v}_B^T \boldsymbol{A}_B \boldsymbol{v}_B \, .
\end{equation}
In this formulation $\boldsymbol{v}_B$ denotes a vector whose length corresponds to the size of block $B$. $\boldsymbol{R}_B^T\boldsymbol{v}_B$ is a prolongated vector whose length corresponds to the size of the full system and which only has nonzero entries at the indices of block $B$. When each index is in at least one block, for every vector $\boldsymbol{v}$ there exists a set $\{\boldsymbol{v}_B\}$ such that $\boldsymbol{v}=\sum_{B \in \mathcal{B}} \boldsymbol{R}_B^T \boldsymbol{v}_B$. With overlapping blocks, different sets $\{\boldsymbol{v}_B\}$ with this property exist. The lemma states that the inner product of a vector $\boldsymbol{v}$ with $\boldsymbol{S}^{-1}$ is equal to the minimum of the sum of inner products over all sets $\{\boldsymbol{v}_B\}$ that sum up to $\boldsymbol{v}$.

Next, recall the mechanism of small eigenvalues in FCM when basis functions on small cut elements are small or almost linearly dependent. The Additive-Schwarz lemma establishes that a vector $\boldsymbol{v}$ that has a small inner product with $\boldsymbol{A}$ because of this mechanism, will also have a small inner product with $\boldsymbol{S}^{-1}$ provided that functions that are almost linearly dependent are aggregated in one block. This clustering can be achieved by assigning one block for every element, consisting of all basis functions that are supported on that element, since basis functions need intersecting supports in order to be almost linearly dependent. When the blocks are set up in this manner, it so follows from the Additive-Schwarz lemma that $\boldsymbol{S}^{-1}$ inherits the small eigenvalues of $\boldsymbol{A}$ caused by small cut elements. The preconditioner can therefore effectively treat the fundamental cause of ill-conditioning of the finite cell method, which has already been demonstrated for uniform discretizations in \cite{CbAS}.

\subsection{Preconditioning the multi-level $hp$-refined finite cell method}\label{sec::preconditioningmlhprefinement}

It is not straightforward to apply this principle to multi-level $hp$-refined grids. To begin with, the notion of an element is less trivial, because of the different levels of elements. Considering the requirement that functions that are almost linearly dependent need to be in the same block, as stated in Section~\ref{sec::additiveschwarzfcm}\,, it is sufficient to only set blocks for the active or leaf elements, i.e.\ elements that do not have overlay elements on a finer level. This still does not yield optimal results however. Therefore the following three modifications are made to the preconditioning technique.

\subsubsection*{Selecting suitable blocks}

Uniform grids, as used in \cite{CbAS}, and multi-level $hp$-grids, as applied here, have two fundamental differences with regard to choosing Additive-Schwarz blocks in an element-wise manner based on the functions supported on a leaf element. First, on uniform grids the number of basis functions supported on an element is constant, i.e.\ $(p+1)^d$ for scalar problems with $p$ denoting the polynomial order and $d$ the number of dimensions \footnote{ For the full tensor product space, e.g\ \cite{Szabo2004}.}\,. For multi-level $hp$-grids the total number of basis functions supported on a leaf element is, in general, not bounded. This is attributed to the superposed linear hat functions that are present on every refinement level. This superposition results in multiply-defined linear basis functions within a leaf element as shown in Figure \ref{fig::multilevelhp::multilevelhpIdea}\,. Directly applying the strategy from \cite{CbAS} to leaf elements in multi-level $hp$-grids therefore yields large blocks that are computationally expensive in the setup, storage and application of the preconditioner. Second, on uniform grids the number of elements that a basis function is supported on is bounded, i.e.\ $\leq2^d$ for $C^0$-finite elements. With multi-level $hp$-discretizations, it is possible that a (coarse) linear basis function is supported on a large number of (fine) leaf elements, Figure \ref{fig::multilevelhp::multilevelhpIdea}\,. Consider the case that a basis function $\phi_k$ is in $n$ different blocks with $n \gg 1$. The unit vector $\boldsymbol{e}_k$ corresponding to $\phi_k$ has a Rayleigh quotient $A_{kk}$ with $\boldsymbol{A}$. Next we construct a set $\{\boldsymbol{v}_B\}$ with the value $1/n$ for all vector entries corresponding to $\phi_k$ and the value $0$ for all other entries. Clearly $\sum_{B \in \mathcal{B}} \boldsymbol{R}_B^T \boldsymbol{v}_B = \boldsymbol{e}_k$ and from Eq.\,\eqref{eq:AdditiveSchwarzLemma} it follows that the Rayleigh quotient of $\boldsymbol{e}_k$ with $\boldsymbol{S}^{-1}$ is bounded from above by $A_{kk}/n$. This shows that a basis function that is in many different blocks yields a small eigenvalue in $\boldsymbol{S}^{-1}$, relative to $\boldsymbol{A}$, and consequently a large eigenvalue in $\boldsymbol{S}$. This reduces the efficiency of the preconditioner, which is demonstrated in Sections~\ref{sec::ballinbox} and \ref{sec::gearbox}\,.

In order to obtain a robust and efficient Additive-Schwarz preconditioner for the multi-level $hp$-finite cell method, the  blocks in $\mathcal{B}$ have to be carefully selected to ensure that \emph{(i)} the number of blocks a given basis function belongs to is bounded and that \emph{(ii)} the size of each block is small. As already noted in \cite{CbAS} it is \emph{sufficient} but not \emph{necessary} to devise a block for every element with \emph{all} basis functions supported on it. In fact it is enough to overcome the conditioning problems of FCM, due to small eigenvalues on small cut elements, if the blocks contain the basis functions that can potentially become almost linearly dependent. This obviously reduces the size of the blocks, and thereby repairs the large computational cost of setting up, storing and applying the preconditioner. Furthermore, the large linear basis functions that are supported on many leaf elements, and cause the second aforementioned efficiency problem, are not the basis functions that suffer from almost linear dependencies. Therefore, truncating the blocks to only the basis functions that can become almost linearly dependent when the element is cut resolves both mentioned issues related to the preconditioning of the multi-level $hp$-finite cell method.

\begin{figure}[H]
 \begin{center}
 \includegraphics[width=0.9\textwidth]{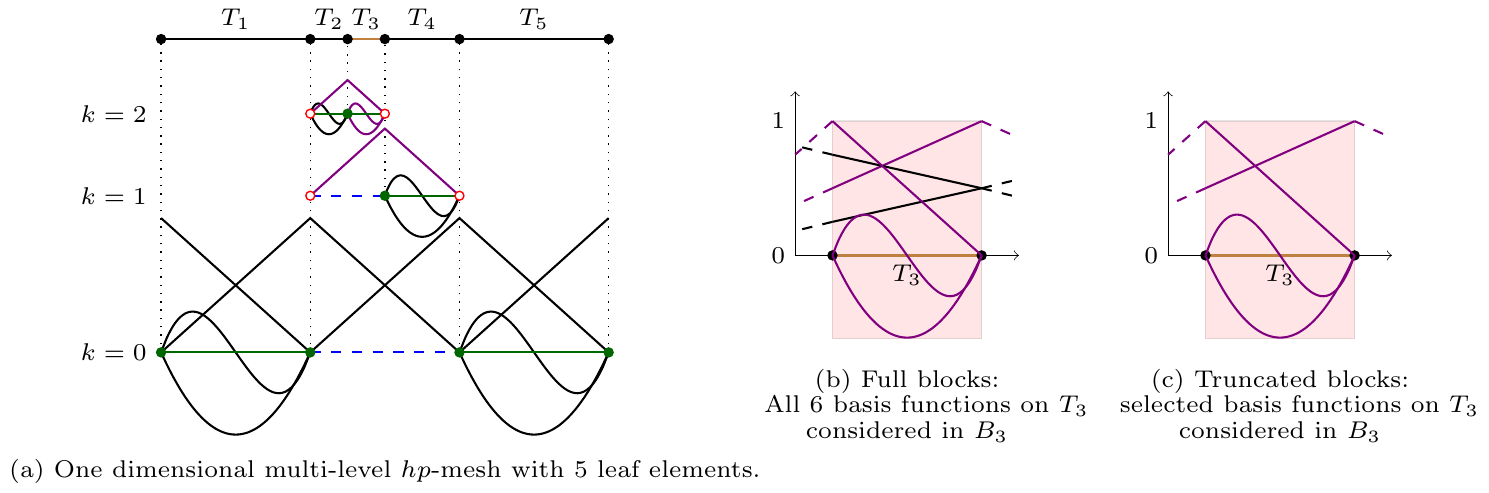}
\caption{Block selection for multi-level $hp$-grids with the full and truncated block of element $T_3$. \label{fig:truncation}}
 \end{center}
\end{figure}

To describe the truncation, we first consider a one dimensional multi-level $hp$-grid with elements of polynomial order $p=3$. Figure~\hyperref[fig:truncation]{\ref*{fig:truncation}a} illustrates that more than $p+1$ basis functions can be supported on a leaf element due to multiply-defined linear basis functions. As mentioned previously, these large sets of supported basis functions need to be truncated in order to result in efficient element-wise Additive-Schwarz blocks. To achieve this, it is first noted that although more then $p+1$ basis functions can be supported on a single element, each element only has $p+1$ unique polynomial degrees of freedom. Since eigenvalues with problematically small eigenmodes are only supported on a small cut element, efficient preconditioning only requires the Additive-Schwarz block of an element to span these $p+1$ polynomial degrees of freedom. This set of $p+1$ functions must be made up of the $p-1$ higher order basis functions accompanied by two linear basis functions. These linear basis functions depend on the position and refinement depth of a leaf element as illustrated in Figure ~\hyperref[fig:truncation]{\ref*{fig:truncation}a}\,. In the grid shown, two different cases are distinguished: \emph{(i)} the element simply has two linear functions on the highest level e.g.\ $T_1$ and $T_5$; \emph{(ii)} the element has one linear function on the highest level and the second linear function is selected by traversing down the element hierarchy and taking the linear function that is linearly independent of the first selected linear function. This can be from a refinement level immediately below the highest level, as is the case for $T_3$ and $T_4$, or from an even lower level in the element hierarchy, as is the case for $T_2$. Figures~\hyperref[fig:truncation]{\ref*{fig:truncation}b} and Figure~\hyperref[fig:truncation]{\ref*{fig:truncation}c} show the block $B_3$ on element $T_3$ before and after truncation, respectively. Note that when the domain in Figure~\hyperref[fig:truncation]{\ref*{fig:truncation}a} is cut either from the left or the right on any of the elements, the truncated Additive-Schwarz block described here always spans all possible functions that are only supported on the cut element and form the problematic eigenmodes for the conditioning of FCM. The truncation in multiple dimensions is achieved by the tensor product of the procedure in one dimension, yielding blocks of $(p+1)^d$ basis functions. Moreover, this choice of blocks not only reduces the size of each block, but also ensures that the number of blocks a single basis function belongs to is at most $2^d$. For clarity it is indicated in all numerical examples in Section~\ref{sec::NumericalResults} whether the full or truncated blocks are applied.

\begin{remark}
Multi-level $hp$-bases are well-conditioned for boundary fitted meshes, due to the orthogonality properties of the superposed linear basis functions and the integrated Legendre basis functions \cite{Zander2016}. In an immersed setting these properties diminish, and when no fictitious stiffness is applied (i.e.\ $\alpha =0$) even full linear dependencies of the linear basis functions can occur in specific cut scenarios. This generates a nullspace resulting in a singular matrix $\mathbf{A}$, and therefore restricts direct solvers to strictly $\alpha > 0$. Iterative solution methods based on Krylov subspaces neglect these nullspaces. The PCG solver that is applied here minimizes the energy in the range $R(\boldsymbol{A})$, e.g.\ \cite{Kaasschieter1988,Saad2003}. Even when $\boldsymbol{A}$ is singular, the right hand side $\boldsymbol{b}$ lies in the range $R(\boldsymbol{A})$, see Eq.\,\eqref{eq::modifiedweakform}. Therefore (preconditioned) iterative solvers are robust to setting $\alpha =0$.
\end{remark}

\subsubsection*{Optimizing the number of blocks}\label{sec::optimization}

The bulk of the computational effort in the construction of the preconditioner consists of inverting the sub-matrices in Eq.\,\eqref{eq:AdditiveSchwarz} and storing the inverses, and is proportional to the number of blocks in $\mathcal{B}$. When a block is devised for \emph{every} element in the mesh, this leads to a large computational cost. Repairing the specific ill-conditioning effects of FCM does not require a block for \emph{every} element, as can be explained by the Additive-Schwarz lemma in Section~\ref{sec::additiveschwarzfcm} and is demonstrated in \cite{CbAS} by only devising blocks for cut elements. In this work we explore whether the computational cost can be further reduced by only devising blocks for cut elements $T$ whose volume fraction $\eta_T$
\begin{equation}
 \eta_T = \frac{|T \cap \Omega_{\text{phys}}|}{|T|} \, , 
\end{equation}
is smaller than a threshold $\bar{\eta}$. Note that this generally results in basis functions in the physical domain that are not contained in any of the blocks, for which Jacobi preconditioning is applied. In this case, the Jacobi preconditioning can technically be interpreted as devising a separate block for every such function. Choosing a smaller value of $\bar{\eta}$ reduces the number of blocks, which yields a more sparse preconditioner and thereby reduces the computational cost. However, smaller values of $\bar{\eta}$ allow smaller untreated elements and with that large condition numbers in accordance with Eq.\,\eqref{eq::kappacut}. This reduces the effectiveness of the preconditioner and increases the required number of iterations in an iterative solver. The effect of $\bar{\eta}$ on the memory usage, number of iterations and computation time is investigated in Section~\ref{sec::etastudybone}\,. 

\subsubsection*{Stabilization using pseudo inverses}\label{sec::stabilityPrec}

The construction of the preconditioner in Eq.\,\eqref{eq:AdditiveSchwarz} requires the inversion of sub-matrices. All these sub-matrices should theoretically be symmetric positive definite, because of the coercive and symmetric weak form in Eq.\,\eqref{eq::modifiedweakform}. Since cut elements can be arbitrarily small, these sub-matrices can contain arbitrarily small eigenvalues, however. When the difference between the largest and smallest eigenmode in a sub-matrix is of the order of the machine precision or smaller, a small eigenvalue that should theoretically be positive can become negative because of rounding errors due to finite precision. Inverting a sub-matrix that contains a negative eigenvalue with an arbitrarily small magnitude yields an inverse with an arbitrary large negative eigenmode. As will be studied numerically in more detail in Section~\ref{sec::vertebra}\,, Conjugate Gradient solvers may diverge when such a mode becomes a dominant part of the error. This problem can be remedied by replacing the inverse of $\boldsymbol{A}_B$ with a \emph{pseudo} inverse. Because of the symmetric positive definite nature of $\boldsymbol{A}_B$ we can write
\begin{equation}
  \boldsymbol{A}_B = \sum_{k=1}^{m_B} \lambda_k \boldsymbol{v}_k \otimes \boldsymbol{v}_k \hspace{3mm} , \hspace{3mm}  \boldsymbol{A}_B^{-1} = \sum_{k=1}^{m_B} \frac{1}{\lambda_k} \boldsymbol{v}_k \otimes \boldsymbol{v}_k \, ,
\end{equation}
with $m_B$ the dimension of $\boldsymbol{A}_B$ and $\lambda_k$ and $\boldsymbol{v}_k$ the $k$-th eigenvalue and eigenvector of $\boldsymbol{A}_B$, respectively. The pseudo inverse that replaces $\boldsymbol{A}_B^{-1}$ and stabilizes the preconditioning technique is defined as
\begin{equation}\label{eq:pseudoinverse}
  \boldsymbol{A}_{B+}^{-1} = \sum_{k=1}^{m_B} \lambda_{k+}^{-1} \boldsymbol{v}_k \otimes \boldsymbol{v}_k \, ,
\end{equation}
with
\begin{equation}
  \lambda_{k+}^{-1} = \begin{cases} \lambda_k^{-1} & \text{for } \lambda_k > \lambda_\text{tres} \, , \\ 0 & \text{for } \lambda_k \leq \lambda_\text{tres} \, , \end{cases}
\end{equation}
and $\lambda_\text{tres} = \varepsilon \max_k (\lambda_k)$ with $\varepsilon$ here selected as $ = 10^{-13}$ such that the smallest eigenvalue that is inverted is sufficiently larger than machine precision.  Although the proposed algorithm requires an eigenvalue decomposition, see Eq.\,\eqref{eq:pseudoinverse}, which is more computationally intensive than the computation of a standard inverse, this routine is not a bottleneck and the setup can be easily accelerated through parallelization as elaborated in Section~\ref{sec::implementation}\,. In Section~\ref{sec::vertebra} the stabilized and unstabilized preconditioner are compared, and it is demonstrated that the stabilized preconditioner does not suffer from the convergence problems of the Conjugate Gradient solver that the unstabilized preconditioner suffers from. All other computations in Section~\ref{sec::vertebra} only apply the stabilized preconditioner.

\begin{remark}
It should be noted that the stabilized preconditioner with pseudo inverses generally contains a nullspace. Therefore iterative solvers based on Krylov subspaces such as the Conjugate Gradient method theoretically do not fully converge to the exact solution, unless $\boldsymbol{b}$ is exactly in the range $R(\mathbf{A})$, e.g.\ \cite{Kaasschieter1988,Saad2003}. However, even if these modes are a part of the solution, this effect is very small in practice because the eigenvalues of the modes disregarded are at least a factor $\varepsilon$ smaller than the largest eigenvalues in the system. To quantify this effect: \emph{(i)} the relative residual can still converge to at least $\varepsilon$ (which is set to be smaller than the tolerance of the solver), \emph{(ii)} the relative energy error as defined in Section~\ref{sec::NumericalResults} still converges to at least $\sqrt{\varepsilon}$ and \emph{(iii)} the relative preconditioned residual can theoretically still fully converge to $0$.
\end{remark}

%% file: algorithm/algorithm.tex
\renewcommand{\rootDir}{algorithm}
\renewcommand{\graphDir}{\rootDir/graphs/}
\renewcommand{\dataDir}{\rootDir/data/}
\renewcommand{\picsDir}{\rootDir/pics}

\section{Implementation and parallelization}\label{sec::implementation}

A major advantage of the proposed preconditioning technique is its suitability for parallel computing, as operations for each block $B$ can be performed independently with minimal synchronization. This can be done using both shared and distributed memory parallelism, as described in the following sections.

\subsection{Shared memory parallelism}\label{sec::openmpparallelism}
Different shared memory paradigms can be used to accelerate the computation of $\boldsymbol{S}$ on shared memory systems. We opt for the use of OpenMP pragmas \cite{openmp08} to speed up the construction by parallelizing the loop over the basis function blocks shown in lines 3 to 7 of Algorithm \ref{alg:preconditioningbadlycutcells}\,. A short study showing the OpenMP scalability of the computation of $\boldsymbol{S}$ will be discussed in Section~\ref{sec::vertebra}\,.

\subsection{Distributed memory parallelism}\label{sec::mpiparallelism}
Distributed memory systems are essential for overcoming the time and memory bottlenecks related to the computation of large engineering problems \cite{Ferronato2012}. In such computations, the computational mesh $\mathcal{T}$ is divided over a set of processes that communicate with one another at well-defined synchronization points via message passing libraries such as the Message Passage Interface (MPI) \cite{mpi1994}. The locality of the developed preconditioner allows it to be easily used in a distributed memory setting. Each MPI-process is identified by its \emph{rank} and assigned a portion $\mathcal{T}_{\textrm{rank}}$ of the computational mesh, with $\mathcal{T}_{\textrm{rank}} \subset \mathcal{T}$. We adopt a hybrid approach combining MPI and OpenMP, where each MPI-process can use OpenMP threads to speed up local computations on every $\mathcal{T}_{\textrm{rank}}$. 

Although iterative solvers are easier to parallelize than direct solvers, e.g.\,\cite{Dongarra:1998:NLA:552704,Saad2003}, they face the challenge of robustness w.r.t.\ the number of processes used \cite{Ferronato2012}. A typical caveat of parallel preconditioned iterative solvers is an increase in the number of iterations needed for convergence when the number of processes utilized in a simulation increases. Because of the sensitivity of Additive-Schwarz preconditioning for the finite cell method caused by the inversion of small cut element contributions, we apply two layers of ghost elements (elements present on a process but owned by another process) to ensure that the exact theoretical preconditioner is obtained regardless of the number of processes or the topology of the partitions $\mathcal{T}_{\textrm{rank}}$. This has the additional advantage that the number of iterations needed for convergence is independent of the number of processes used and thus favors good parallel scalability. It should be noted that the evaluation of the ghost element layers does not significantly increase the setup time of the linear system, since the number of ghost elements on a process is much smaller than the number of active elements. Moreover, since the bulk of the execution time is spent in solving the linear system, the iterative solver's scalability largely governs the overall scalability of the simulation.

\begin{figure}[H]
   \begin{center}
             \includegraphics[width=0.8\textwidth]{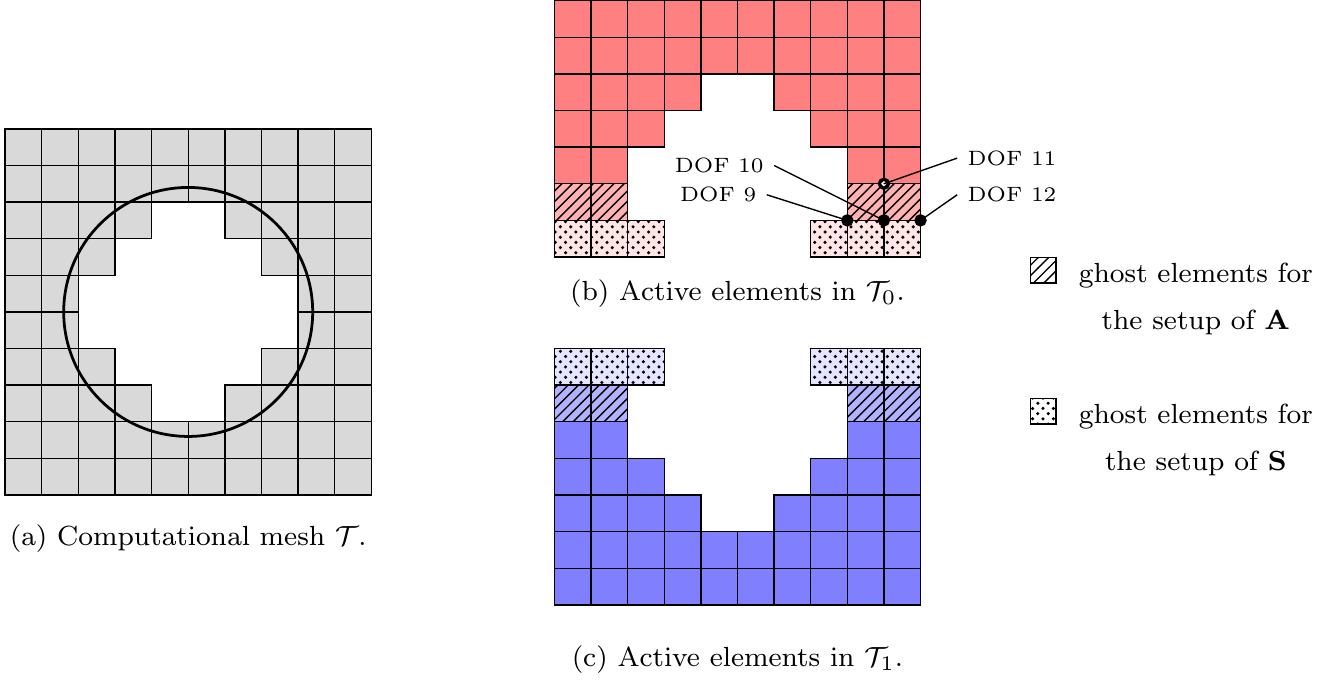}
   \hfill%
        \caption{Partitioning of the computational mesh in parallel computations for the setup of $\boldsymbol{A}$ and $\boldsymbol{S}$.}
        \label{fig::mpiparallelism}
    \end{center}
\end{figure}

Figure \ref{fig::mpiparallelism} depicts the partitioning of a computational mesh between two MPI-processes. The red and blue elements represent active elements on $\mathcal{T}_0$ and $\mathcal{T}_1$ owned by process 0 and process 1 respectively. The hatched elements constitute the first layer of ghost elements. This layer allows a process to compute the full contributions in $\mathbf{A}$ associated with DOFs it owns without communicating with other processes. This is illustrated in Figure~\hyperref[fig::mpiparallelism]{\ref*{fig::mpiparallelism}b} where process 0 can independently compute all entries in $\mathbf{A}$ associated with DOF 11 by integrating all the elements around this node.   
Likewise, the second layer of ghost elements ensures that the exact theoretical preconditioner $\boldsymbol{S}$ can be constructed without communication. Figure~\hyperref[fig::mpiparallelism]{\ref*{fig::mpiparallelism}b} also illustrates the importance of this second ghost-element layer. Since DOF 11 shares a support with DOFs 9, 10 and 12, process 0 requires the complete entries in $\mathbf{A}$ associated with these DOFs to ensure that all entries in $\mathbf{S}$ related to DOF 11 are computed independent of the partitioning as per Eq.\,\eqref{eq:AdditiveSchwarz}. This guarantees that the convergence of an iterative solution method is independent of the number of processes and how the elements are divided over these.
The above explanation considers uniform meshes, but can trivially be extended to locally refined grids. The parallel performance of a preconditioned parallel Conjugate Gradient solver using the described approach will be investigated in Section \ref{sec::mpiscalability}\,. 

\subsection{Algorithm}

The preconditioner construction is summarized in Algorithm \ref{alg:preconditioningbadlycutcells}\,. The volume fraction $\eta_T$ of an element is used as a criterion to decide whether a block $B$ is devised for the basis functions supported on it. Note that the algorithm does not explicitly devise a separate block of one basis function for the basis functions that are not present in any of the truncated blocks $B \in \mathcal{B}$ --- as the diagonal scaling of these basis functions can be interpreted theoretically, Section~\ref{sec::preconGeneral} --- but that a separate diagonal scaling step is done for these basis functions in lines $8$ to $14$.

\begin{algorithm}[H]
	\caption{Computation the preconditioner $\boldsymbol{S}$}
\label{alg:preconditioningbadlycutcells}
\begin{algorithmic}[1]
 \Function{\texttt{computePreconditioner} }{ $\boldsymbol{A}$, $\{\eta_T\}$, $\bar{\eta}$ }
  \State  $\mathcal{B}$  = getBlocks( $\boldsymbol{A}$, $\{\eta_T\}$, $\bar{\eta}$ )
  \State  \textbf{for} $B \in \mathcal{B}$ \textbf{do} \Comment( Loop over basis function blocks )  
	 \State \hspace{3mm} $\boldsymbol{A}_{B}$ =  extractSubMatrix( $\mathbf{A}$, $B$ )
	 \State \hspace{3mm} $\boldsymbol{A}_{B+}^{-1}$ =  computePseudoInverse( $\boldsymbol{A}_{B}$ ) 
  \State \hspace{3mm} $\boldsymbol{S}$  +=  prolongateMatrix( $\boldsymbol{S}, \boldsymbol{A}_{B+}^{-1}$) 
  \State \textbf{end} 
  \State  $n$ = size$(\boldsymbol{A})$
  \State  $l=1$
  \State \textbf{for} $l < n+1$ \textbf{do} \Comment( Diagonal scaling ) 
  \State \hspace{3mm} \textbf{if} ( $\boldsymbol{S}(l,l) == 0$ ) \textbf{do} 
  \State \hspace{6mm} $\boldsymbol{S}(l,l)$ = $\dfrac{1.0}{\boldsymbol{A}(l,l)}$ 
  \State \hspace{3mm} $l$ += $1$
  \State \textbf{end} 
  \State \hspace{-5mm}\textbf{return} $\boldsymbol{S}$
 \EndFunction
\end{algorithmic}
\end{algorithm}

%% file: numericalresults/numericalresults.tex
\clearpage
\section{Numerical examples}\label{sec::NumericalResults}

\renewcommand{\rootDir}{numericalresults}
\renewcommand{\graphDir}{\rootDir/graphs}
\renewcommand{\dataDir}{\rootDir/data}
\renewcommand{\picsDir}{\rootDir/pics}

This section demonstrates the suitability of the proposed preconditioner for finite cell analysis of real-life problems. The numerical examples are specifically chosen to highlight various aspects presented in the previous sections. Focus is placed on image-based geometries (in Sections \ref{sec::vertebra} and \ref{sec::gearbox}), an application field where the finite cell method is particularly advantageous  in circumventing laborious mesh generation procedures. The first example addresses the extension of the Additive-Schwarz preconditioning technique to FCM problems involving multi-level $hp$-refinement by studying a simple example of a cube with a spherical cavity subjected to uniaxial loading. The second example, from the field of biomechanics, considers the loading of a lumbar vertebral body and studies the effect of the pseudo inverse and the volume fraction threshold $\bar{\eta}$ on the computational cost of a preconditioned Conjugate Gradient solver (PCG). The final example brings together FCM, multi-level $hp$-refinement, Additive-Schwarz preconditioning and parallel computing in the image-based analysis of a die cast gearbox housing.

The convergence of the PCG solver applied in this section is assessed by monitoring the development of the relative preconditioned residual
\begin{equation}\label{eq:residual}
 \frac{\|\mathbf{S}  \boldsymbol{r}_i\|}{\| \mathbf{S} \boldsymbol{b}\|} = \frac{\|\mathbf{S} \boldsymbol{b} - \mathbf{S} \boldsymbol{A} \boldsymbol{x}_i\|}{\| \mathbf{S} \boldsymbol{b}\|} \, ,
\end{equation}
and when possible the relative error in the energy norm
\begin{equation}
\label{eq::energyerror}
e_i = \frac{\| \boldsymbol{x}_\text{ref} - \boldsymbol{x}_i \|_{\boldsymbol{A}}}{\| \boldsymbol{x}_\text{ref} \|_{\boldsymbol{A}}} = \frac{\sqrt{\left( \boldsymbol{x}_\text{ref} - \boldsymbol{x}_i \right) \boldsymbol{A} \left( \boldsymbol{x}_\text{ref} - \boldsymbol{x}_i \right)^T}}{\sqrt{\boldsymbol{x}_\text{ref} \boldsymbol{A} \boldsymbol{x}_\text{ref}^T}} \, ,
\end{equation}
with a reference solution $\boldsymbol{x}_\text{ref}$ obtained from the parallel direct solver \texttt{Intel}\textsuperscript{\textregistered} \texttt{Pardiso}, contained in the Intel Math Kernel Library \cite{intel}. The subscript $i$ in Equations \eqref{eq:residual} and \eqref{eq::energyerror} denotes quantities in the $i$-{th} Conjugate Gradient iteration. In \cite{CbAS} the effectiveness of the preconditioning technique is assessed by focusing on the condition number of the system. We do not adopt this approach here but rather focus of the residual and energy convergence since the cost of computing condition numbers for the large systems under investigation is prohibitive.    

\subsection{Compression of a cube with a spherical exclusion}\label{sec::ballinbox}

The following example illustrates the importance of selecting suitable blocks for Additive-Schwarz preconditioning in FCM simulations with multi-level $hp$-refinement as described in Section \ref{sec::preconditioningmlhprefinement}\,. To this end, a simple example comprised of a cube with a spherical cavity under compressive loading is considered. Although this setup yields a relatively small system, with a small number of DOFs ($< 17\,000$), it serves as a good starting point to illustrate the effectiveness of the preconditioning technique developed in this manuscript as it already shows the conditioning problems associated with FCM and multi-level $hp$-refinement. 

\begin{figure}[H]
 \begin{center}
  \includegraphics[width=0.38\textwidth]{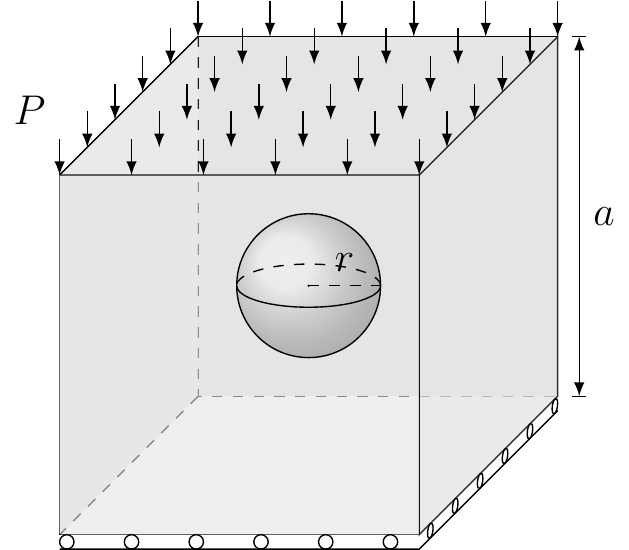}%
  \caption{Cube with a spherical cavity under compressive loading.}
 \label{fig:ballinbox}
 \end{center}
\end{figure}

\subsubsection*{Problem setup}
A cube of unit length with a spherical cavity of radius $r=0.01$ is subjected to a homogeneous pressure load $P$ as shown in Figure \ref{fig:ballinbox}\,. The cube has a Young's modulus of $70$ GPa and Poisson's ratio $\nu=0.34$. Homogeneous Dirichlet boundary conditions in the normal direction are applied using a penalty parameter $\beta=10^{10}$. The base mesh consists of $8 \times 8 \times 8$ elements and is refined toward the surface of the exclusion with a refinement depth $k \in \lbrace  0 , 1, 2, 3, 4 \rbrace$. A polynomial order $p=3$ is used, resulting in the numbers of unknowns for the different refinement depths given in Table \ref{tab:cubedofs}\,. Numerical integration is performed using octree partitioning with a tree-depth of 3 relative to the local element size, see e.g.\ \cite{schillinger_small_2012} for details.   

\begin{table}[H]
\centering
\begin{tabular}{cccccc}
\toprule
$k$ &  0  & 1 & 2 & 3 & 4 \\
\toprule
DOFs & 13\,851 & 14\,541  & 15\,231  & 15\,921  & 16\,611 \\
\bottomrule
\end{tabular}
\caption{Number of degrees of freedom for different refinement depths $k$. }
\label{tab:cubedofs}
\end{table}

\subsubsection*{Convergence behavior}
\begin{figure}[H]
   \begin{center}
	   \subfloat[Full preconditioner: residual convergence.]
     {
       \includegraphics[width=0.48\textwidth]{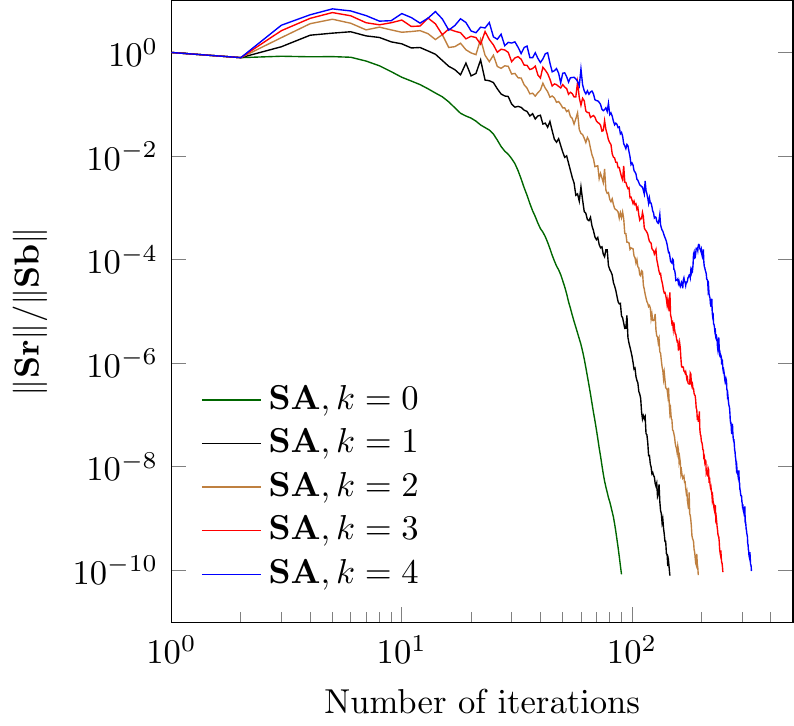}%
       \label{fig:fulleascubeconvergence}
     }%
     \hfill%
	   \subfloat[ Truncated preconditioner: residual convergence. ]
     {
       \includegraphics[width=0.48\textwidth]{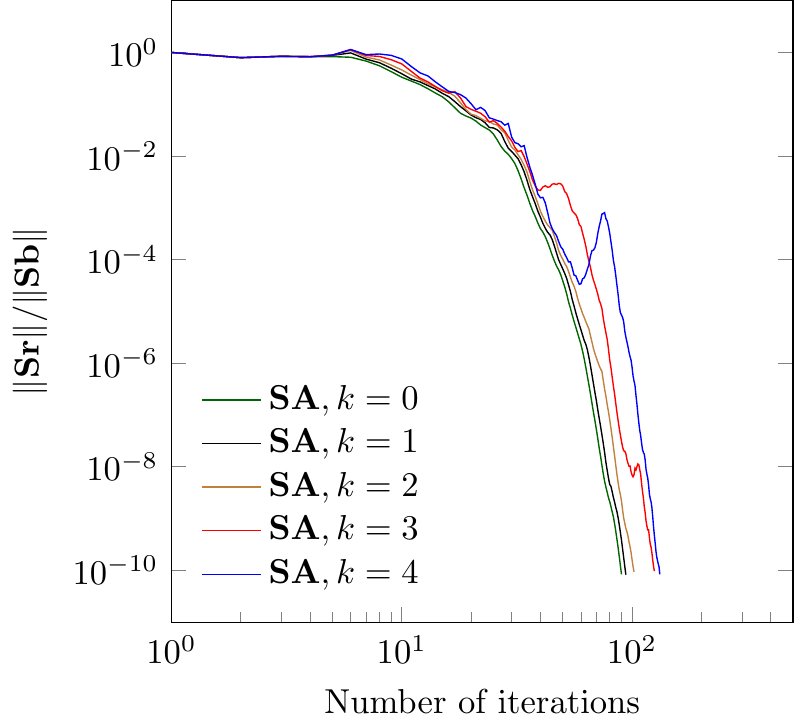}%
       \label{fig:trunceascubeconvergence}
     }%
     \hfill%
	   \subfloat[Full preconditioner: energy convergence.]
     {
       \includegraphics[width=0.48\textwidth]{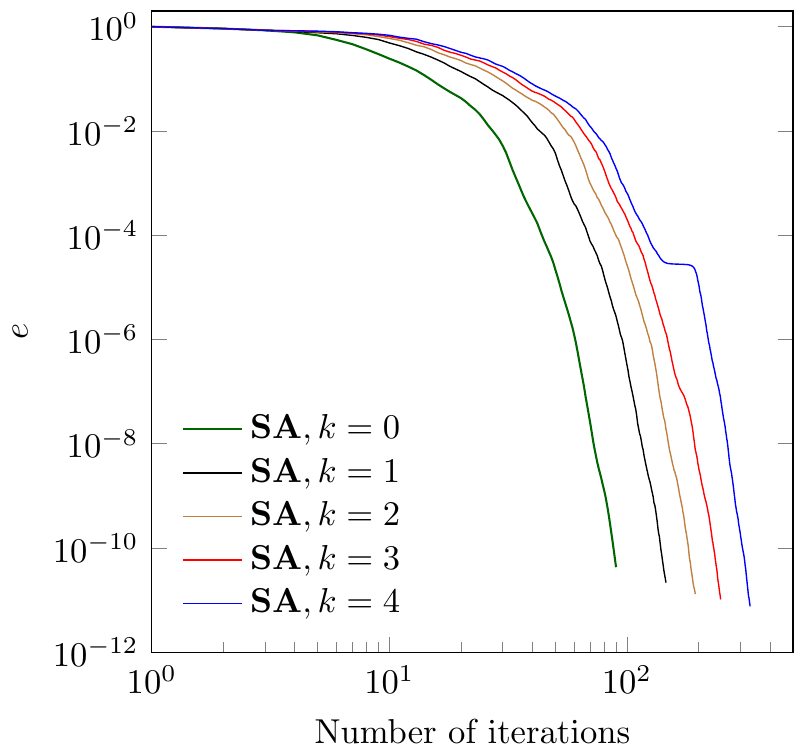}%
       \label{fig:fulleascubeenergy}
     }%
     \hfill%
	   \subfloat[Truncated preconditioner: energy convergence.]
     {
       \includegraphics[width=0.48\textwidth]{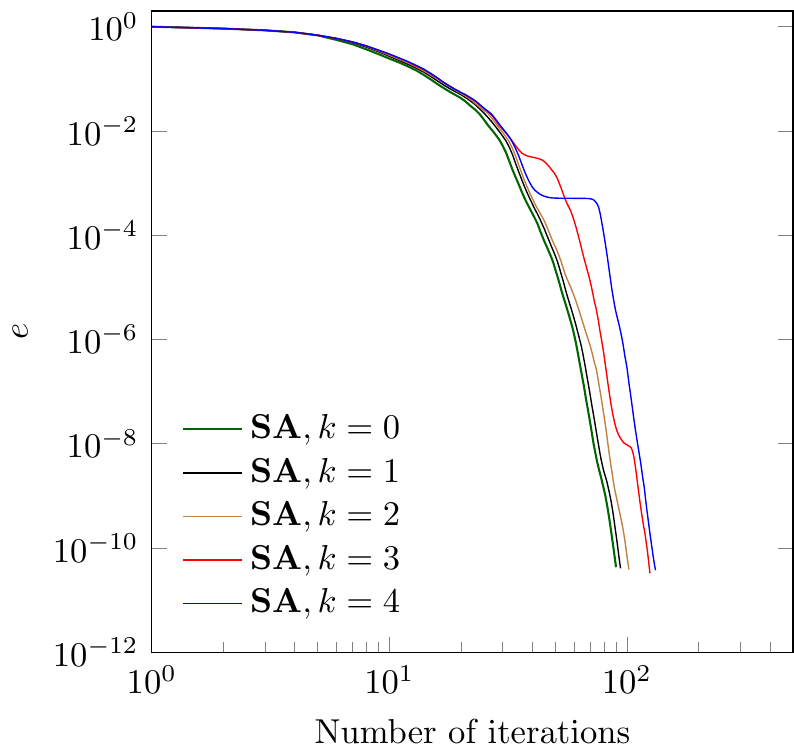}%
       \label{fig:trunceascubeenergy}
     }%
	\caption{ Convergence behavior of the PCG solver with full and truncated blocks for multi-level $hp$-refinement. Note that the convergence of the residual is not monotone as the PCG solver minimizes the energy and the residual is only loosely bound to this, e.g.\ \cite{Saad2003}\,.}
    \label{fig:cubeconvergence}
    \end{center}
\end{figure}

The convergence of the PCG solver preconditioned with the full and the truncated blocks is shown in Figure \ref{fig:cubeconvergence}\,. The results show that the convergence of the residual closely follows the convergence of the energy. This is generally not the case for ill-conditioned systems such as unpreconditioned FCM, e.g.\ \cite{SIPIC}, and indicates that both preconditioners are robust w.r.t.\ small cut cells. Figures \ref{fig:fulleascubeconvergence} and \ref{fig:fulleascubeenergy}\,, however, show the effect of local refinements on the convergence speed of the solver when no truncation is used in the construction of $\mathbf{S}$. This indicates that the preconditioner with full blocks is not robust w.r.t.\ the refinement depth --- i.e.\ a small increase in the number of unknowns, by approximately a factor of 1.2, leads to a significant increase in the number of iterations needed to reach convergence, which increases by a factor of around 3.6. The intensity of this effect is significantly reduced when the preconditioner with truncated blocks is used as demonstrated in Figures \ref{fig:trunceascubeconvergence} and \ref{fig:trunceascubeenergy}\,. The truncated preconditioner is therefore not only computationally less expensive than the full preconditioner (due to the smaller block sizes), but also yields better conditioning. Similar results are obtained using a more complex geometry in the third example (Section \ref{sec::gearbox}).

\subsection{Image-based simulation of a lumbar vertebra}\label{sec::vertebra}

The applicability of the preconditioner in biomechanical simulations is highlighted in the following image-based simulation of a lumbar vertebra subjected to compressive loading, with the computational setup following \cite{Elhaddad2017}. Elastostatic computations are performed with the described setup, utilizing a PCG solver for the solution of the linear system of equations with either no preconditioner, diagonal scaling, or the presented Additive-Schwarz preconditioner. Various aspects of the iterative solver's convergence behavior are monitored during the simulation for different solver and preconditioner configurations. This example does not consider local mesh refinement, which yields a relatively simple problem on a complex geometry, which is a good test case to benchmark the stability of the preconditioning technique and the efficiency of different threshold volume fractions $\bar{\eta}$, Section~\ref{sec::preconditioningmlhprefinement}\,. Note that on the uniform mesh without local refinements truncation of the blocks is not applicable.

\subsubsection*{Problem setup}\label{sec::lumbarsetup}
\begin{figure}[t]
 \begin{center}
  \begin{minipage}[c]{0.5\textwidth}
   \begin{center}
    \includegraphics[width=1.0\textwidth]{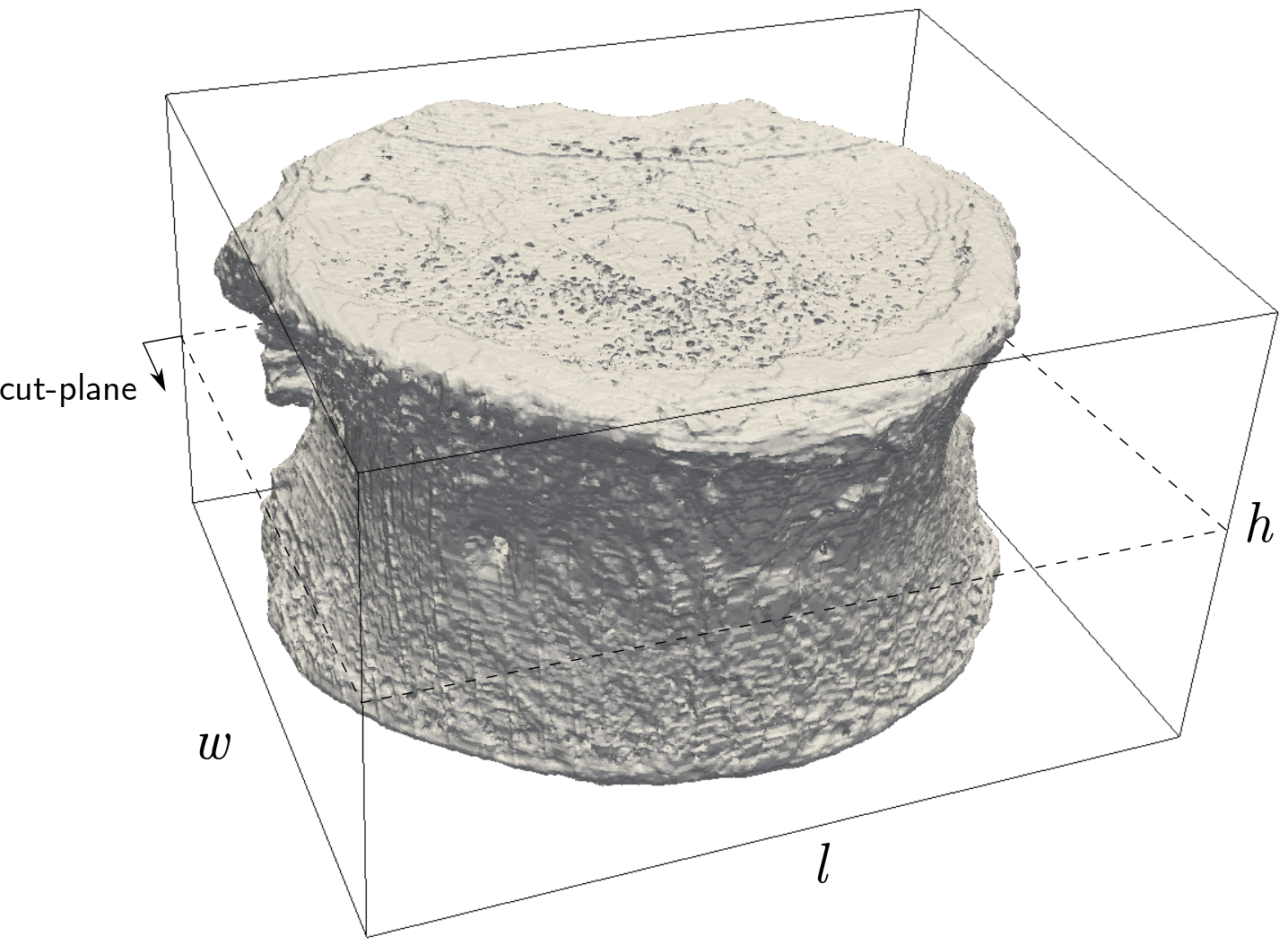}%
   \begin{table}[H]
   \centering
   \scriptsize
   \begin{tabular}{cc}
   \toprule
   Material &  vertebral bone \\
   Young's modulus, $E_{\text{bone}}$& 12 GPa \\
   Poisson's ratio, $\nu$ & 0.3 \\
   $l \times w \times h$ [mm] & $ 47.3 \times 55.3 \times 34.1 $   \\
   \bottomrule
   \end{tabular}
   \end{table}
   
   (a) Geometric and material properties of the lumbar vertebra. 

   \end{center}
  \end{minipage}
 \hspace{2mm}
 \begin{minipage}[c]{0.4\textwidth}
  \scriptsize
  \begin{center}
   \hspace*{1cm}\includegraphics[width=0.8\textwidth]{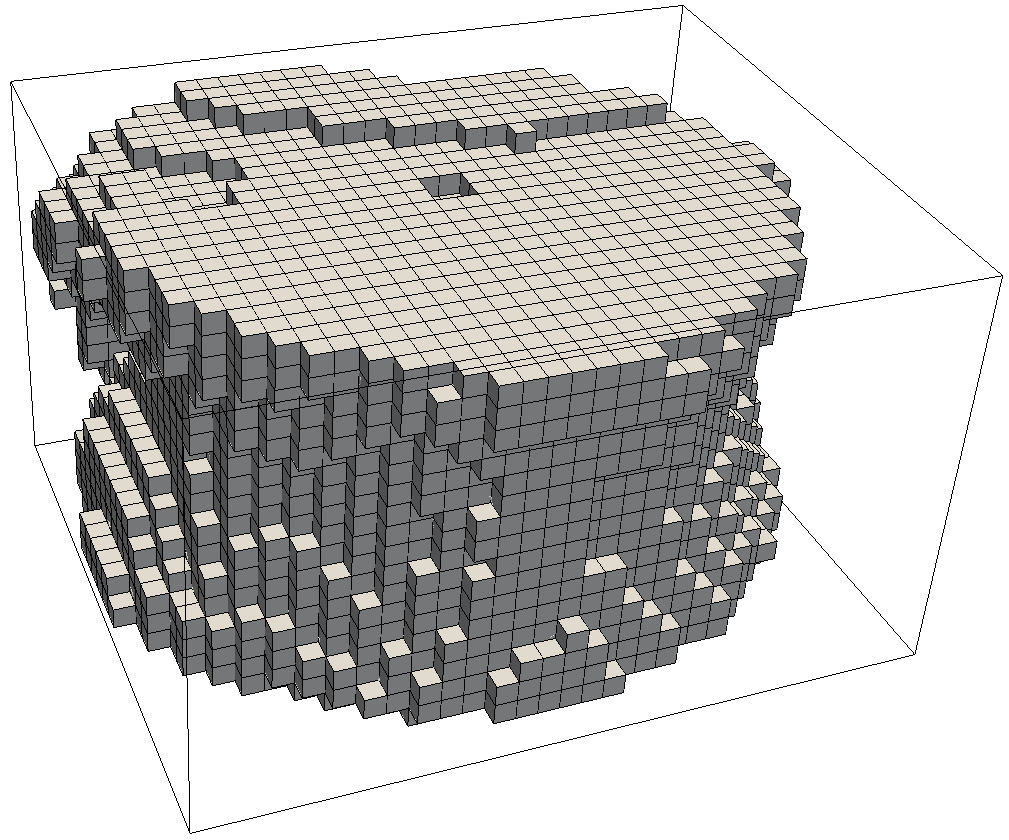}%
   \\ \hspace{0.5cm} (b) Computational mesh of active elements.
   \\
   \vspace{3mm}
   \includegraphics[width=0.7\textwidth]{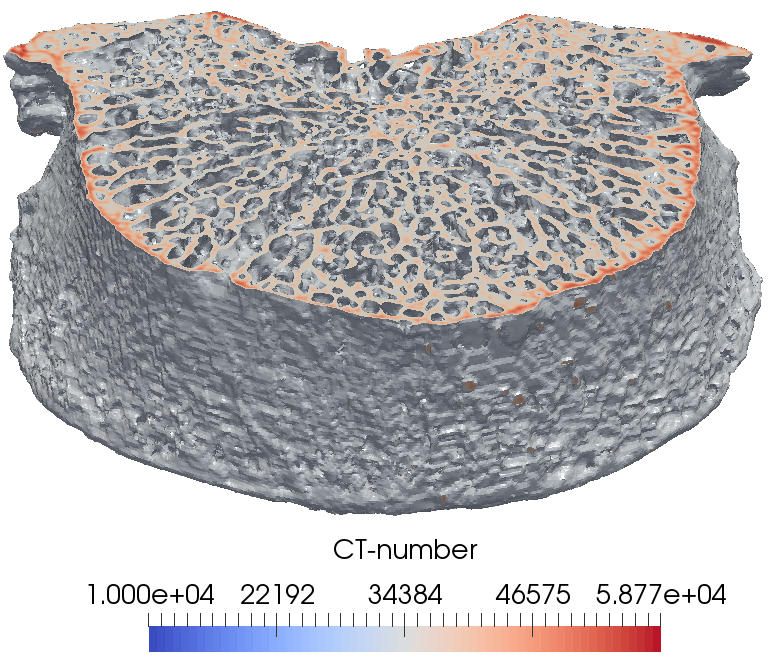}%
   \\ (c) Cut showing the porous bone structure and distribution of the CT-number. 
  \end{center}
  \end{minipage}
  \caption{Setup of the lumbar vertebra example.}
  \label{fig::lumbarvertebra}
 \end{center}
\end{figure}

The geometry of the lumbar vertebra is obtained via a high-resolution micro-CT scan of the specimen with a voxel size of  $80 \times 80 \times 80 \ \mu m^3$. A subsequent segmentation of the scan using \texttt{ITK-SNAP} \cite{py06nimg} yields a model of the vertebral body without the surrounding soft tissue, Figure \hyperref[fig::lumbarvertebra]{\ref*{fig::lumbarvertebra}a}\,. It is noted that the segmentation only considers the outer boundary of the vertebra and no distinction is made between cortical and trabecular bone.  

A non-boundary conforming discretization is generated using the finite cell method with $10^3$ voxels contained in each element resulting in 170\,982 elements. Elements completely outside the physical domain are removed in a preprocessing step yielding 75\,821 active elements. A polynomial order of $p=3$ is chosen, resulting in approximately 1.7 million degrees of freedom. Material properties within each element are defined on a voxel-level by applying a threshold on a voxel's CT-number. The CT-number is a measure of the x-ray absorption coefficient that can be related to the density and consequently the stiffness of the material in a voxel. Material identified as inside the bone structure is assigned a Young's modulus $E_{\text{bone}} = 12$ GPa and Poisson's ratio $\nu=0.3$, corresponding to values commonly used in literature \cite{Keaveny}. Fictitious material, on the other hand, is assigned the material properties $E_{\text{fict}}=10^{-11}$ GPa and $\nu=0.3$ in order to avoid peak stresses in cut elements \cite{NitscheStability}. The surface of the entire vertebral body is recovered using a marching cubes algorithm \cite{Maple} and thereafter trimmed to provide a triangulation for the superior end-plate, where a uniformly distributed vertical load of 800N is applied, and the inferior end-plate, which is weakly clamped using the penalty method with $\beta=10^6$.

\subsubsection*{Convergence behavior and stability of $\mathbf{S}$}
\begin{figure}[t]
   \begin{center}
	   \subfloat[Convergence of the relative residual]
     {
       \label{fig::boneconvergence}
       \includegraphics[width=0.4901\textwidth]{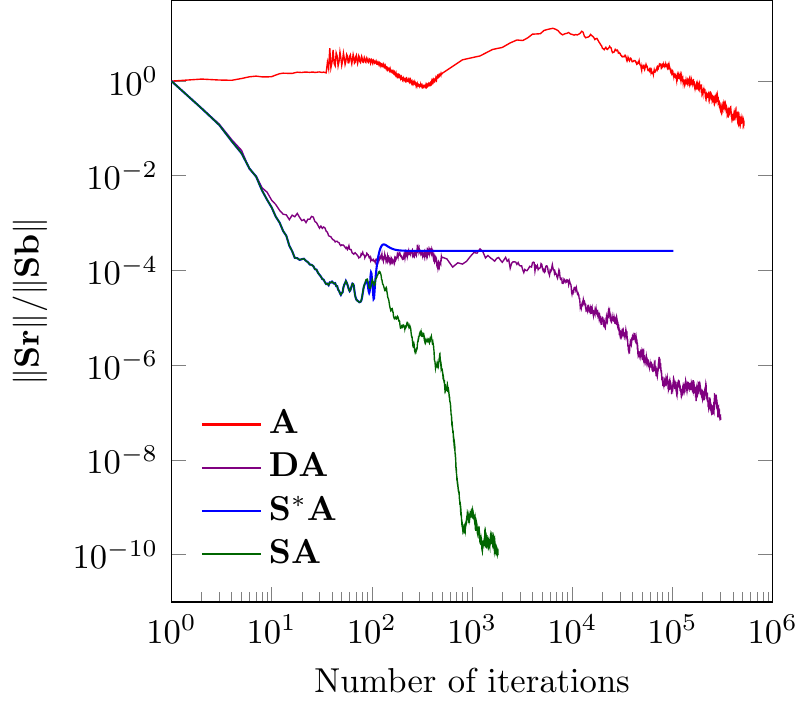}%
     }%
	\hfill%
     \subfloat[Convergence of the error in the energy norm]
     {
       \label{fig::boneenergy}
       \includegraphics[width=0.48\textwidth]{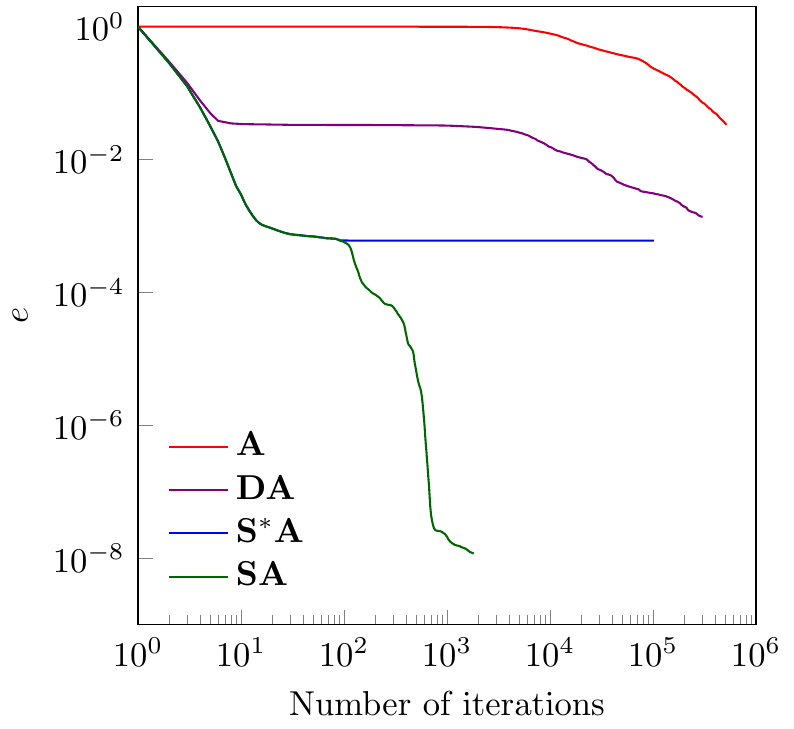}%
     }%
	\caption{ Convergence of the relative residual and the energy norm errors for different preconditioners in the lumbar vertebra example.}
       \label{fig::boneexample}
    \end{center}
\end{figure}

 Figure \ref{fig::boneexample} shows the values of the relative preconditioned residual and the energy error plotted against the number of PCG iterations performed. The solver is terminated when either a cut-off tolerance in the relative preconditioned residual of $10^{-10}$ or the maximum number of iterations (300\,000) is reached. One can observe that the unpreconditioned system, denoted by the curves labeled $\mathbf{A}$, suffers from ill-conditioning characterized by the slow residual convergence rates and high energy norm errors even after a significant number of iterations. The diagonally preconditioned system, $\mathbf{D}\mathbf{A}$, shows an improved convergence behavior but also fails to reach the desired cut-off tolerance within the specified maximum number of iterations.
 
The curves $\boldsymbol{S}^*\boldsymbol{A}$ and $\boldsymbol{S}\boldsymbol{A}$ depict the systems preconditioned with the preconditioning scheme developed here using the full inverse ($\boldsymbol{S}^*$) and stabilized pseudo inverse ($\boldsymbol{S}$) respectively, with cut-off tolerance $\varepsilon=10^{-13}$. The importance of filtering out eigenvalues smaller than the machine precision in the local inverses $\mathbf{A}_{B+}^{-1}$, as described in Section \ref{sec::preconditioningmlhprefinement}\,, is clearly visible from these curves in Figures~\ref{fig::boneconvergence} and \ref{fig::boneenergy}\,. Both $\boldsymbol{S}^*\boldsymbol{A}$ and $\boldsymbol{S}\boldsymbol{A}$ exhibit the same convergence until the eigenvectors of $\boldsymbol{S}^*\boldsymbol{A}$ with negative eigenvalues become the dominant part of the error. It is clearly visible that the relative energy norm error stagnates and the relative residual diverges. These problematic eigenvalues are discarded by $\boldsymbol{S}\boldsymbol{A}$, resulting in convergence up to the desired tolerance as indicated in the figures. It should be noted that potential loss of accuracy due to the discarding of extremely small eigenmodes is not observed in this example.

It is customary in FCM to preprocess the system and discard all functions that completely lie outside of the physical domain. By stabilizing with the pseudo inverse and setting the fictitious stiffness $\alpha = 0$, this preprocessing step is theoretically not required since the pseudo inverse will ignore these functions. Discarding these functions in a preprocessing step still reduces the computational cost however. It is noted that with $\alpha$ large enough the system does not contain eigenvalues of the order of the machine precision, and the pseudo inverse is equivalent to the full inverse. 

\subsubsection*{Optimization: Preconditioning only severely cut cells}\label{sec::etastudybone}

In \cite{CbAS} it is mentioned that the number of blocks used in the setup of the preconditioner can be reduced by considering only elements below a certain volume fraction. This idea is followed in Section \ref{sec::optimization} where it is proposed to only invert blocks of basis functions for elements with a volume fraction smaller than a threshold value $\bar{\eta}\in[ 0.0, 1.0 ]$. Basis functions that are not present in any of these elements are instead diagonally scaled, to reduce the computational cost.

In this section we study the effect of the threshold value $\bar{\eta}$ on the performance of the preconditioner for the lumbar vertebra model described above. To this end, the linear system is solved by a PCG solver for different values of $\bar{\eta}$. A value $\bar{\eta} = 1.0$ corresponds to the preconditioner with a block for every element while $\bar{\eta}=0.0$ results in a Jacobi preconditioner. It is again stipulated that in this computation the mesh does not contain local refinements, such that the blocks are not truncated as described in Section~\ref{sec::preconditioningmlhprefinement}.

\begin{figure}[H]
   \begin{center}
    	   \subfloat[Convergence of the relative residual for various $\bar{\eta}$]
     {
       \includegraphics[width=0.455\textwidth]{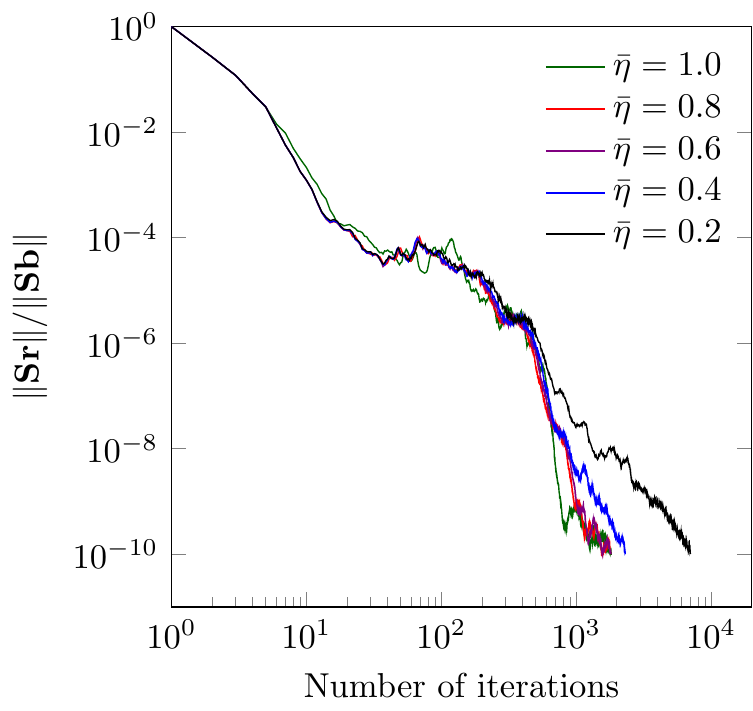}%
       \label{fig::cbasetastudy::residual}
     }%
	\hfill%
	   \subfloat[Convergence of the energy error for various $\bar{\eta}$]
     {
       \includegraphics[width=0.45\textwidth]{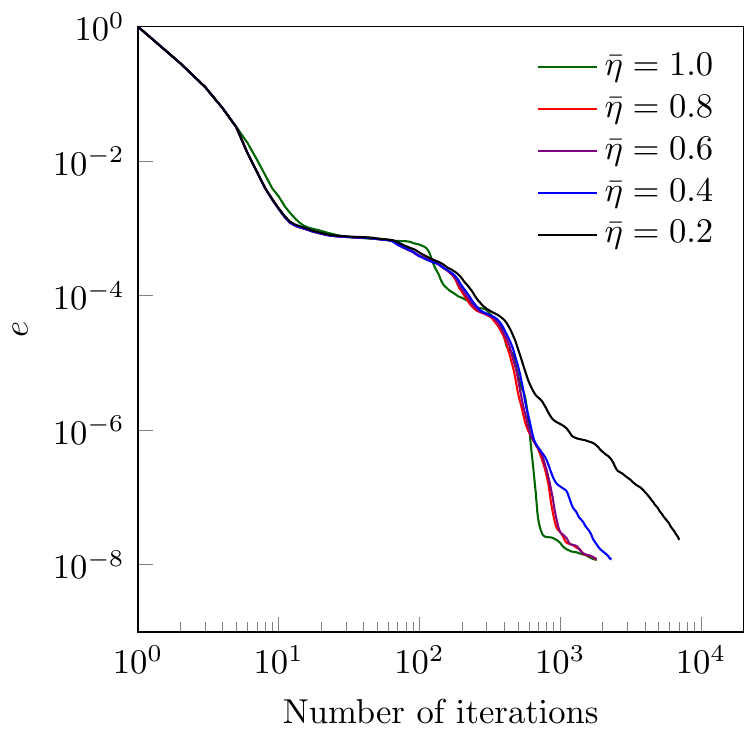}%
       \label{fig::cbasetastudy::energy}
     }%
     \hfill%
     \subfloat[Preconditioner size for different values of $\bar{\eta}$]
     {
       \includegraphics[width=0.46\textwidth]{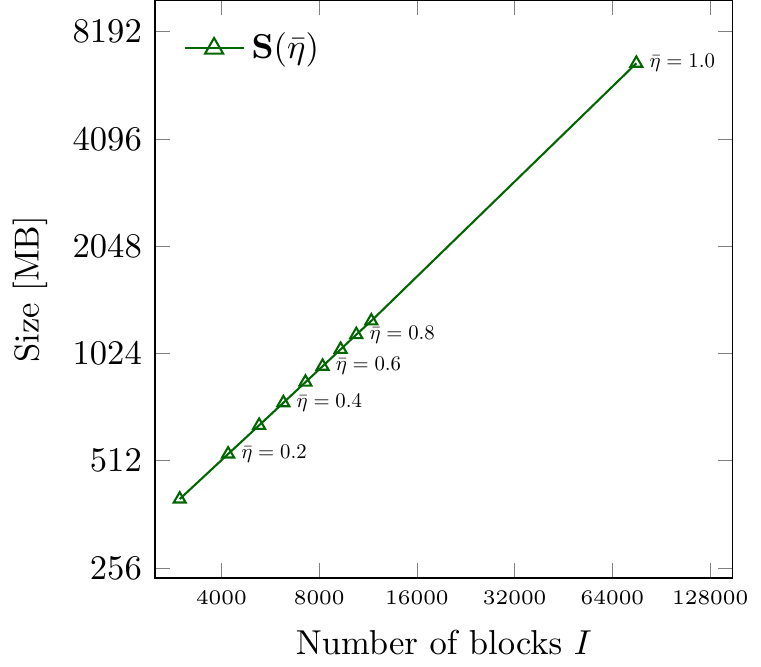}%
       \label{fig::cbasetastudy::sizes}
     }%
     \hfill
	   \subfloat[Convergence behavior for different values of $\bar{\eta}$]
     {
       \put(45,202){$\cdot 10^3$}
       \put(190,202){$\cdot 10^2$}
	   \includegraphics[width=0.51\textwidth]{\graphDir/memoryLSD.tikz}%
       \label{fig::cbasetastudy::convergence}
     }%
	\caption{ Preconditioning of cut cells with $\eta_{T} \leq \bar{\eta}$.}
       \label{fig::cbasetastudy}
    \end{center}
\end{figure}

Figure \ref{fig::cbasetastudy} depicts the influence of $\bar{\eta}$ on the convergence of the PCG solver, the solution time and memory requirements of setting up the preconditioner. Figures \ref{fig::cbasetastudy::residual} - \ref{fig::cbasetastudy::energy} confirm that the effectiveness of the preconditioner, as described in Section \ref{sec::optimization}\,, is maintained as long as problematic basis functions are taken into account when constructing $\mathbf{S}$. This holds for $\bar{\eta} \in [ 0.4, 1.0]$, while values of $\bar{\eta} \le 0.4$ allow smaller untreated elements and thus reduce the effectiveness of $\mathbf{S}$. Reducing the value of $\bar{\eta}$ leads to a sparser preconditioner which requires less storage space as shown in Figure \ref{fig::cbasetastudy::sizes}\,, where the size of the preconditioner scales linearly with the number of blocks in $\mathcal{B}$. The increased sparsity of $\mathbf{S}$ for values of $\bar{\eta} < 1.0$ leads to a reduction in the computational cost of a single PCG iteration, since the matrix-vector multiplication involving $\mathbf{S}$ and $\boldsymbol{r}$ takes less time. This leads to an overall faster execution time as shown in Figure \ref{fig::cbasetastudy::convergence}\,, where the number of iterations needed until convergence is plotted alongside the corresponding solution time for different values of $\bar{\eta}$. The solver converges in less time for $0.6 \leq \bar{\eta} < 1.0$, with the solver converging 1.6 times faster for $\bar{\eta}=0.7$ than for $\bar{\eta}=1.0$. This study shows that values of $\bar{\eta} \in [0.6, 0.7]$ yield an effective preconditioner which is sparser and consequently faster than the \emph{complete} preconditioner for $\bar{\eta}=1.0$ with a block for \emph{every} element. 

To enable comparison of all presented examples, in the remainder of this manuscripts the non-optimized choice of $\bar{\eta} = 1.0$ will be considered. 
\subsubsection*{Computational cost of setting up $\mathbf{S}$} \label{sec::openmpstudy}
\begin{figure}[H]
 \begin{center}
  \includegraphics[width=0.48\textwidth]{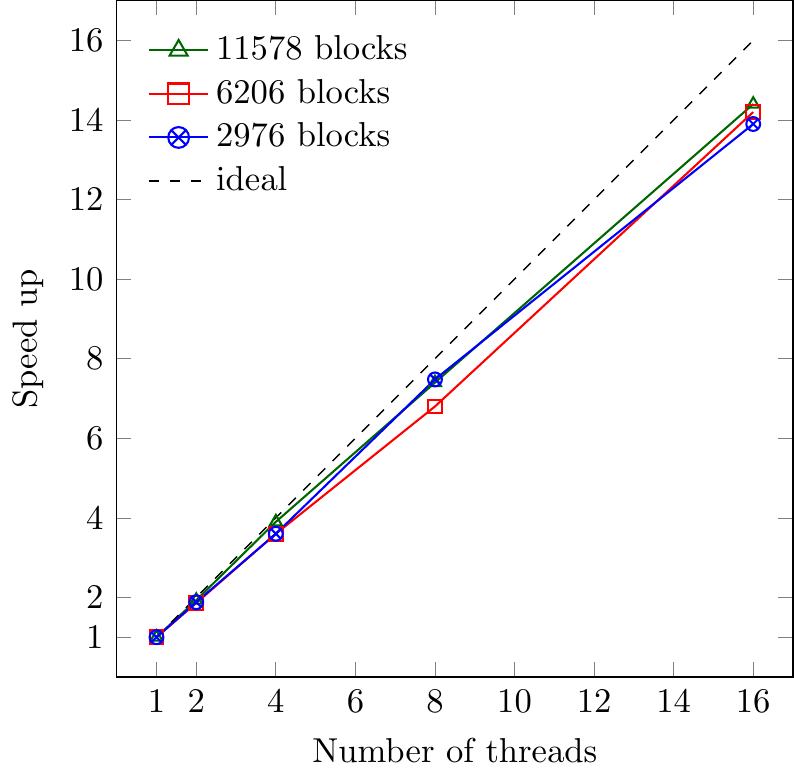}%
  \caption{Setup of $\mathbf{S}$ --- OpenMP scalability for different numbers of blocks in $\mathcal{B}$.}
  \label{fig::ompstudy}
 \end{center}
\end{figure}

Section \ref{sec::openmpparallelism} explains how the construction of $\mathbf{S}$ can be accelerated by the use of shared memory parallelism. This is illustrated in a short study that shows the OpenMP speed up achieved for a variable number of threads involved in the setup of the preconditioner. Three different preconditioners with different numbers of blocks in $\mathcal{B}$, and consequently different computational costs, are analyzed. The results of this study are shown in Figure \ref{fig::ompstudy}\,. Although good parallel scalability is achieved for the different number of blocks, the parallel efficiency decreases gradually with an increase of the number of threads. This is attributed to the synchronization needed when assembling the entries of $\mathbf{A}_{B+}^{-1}$ into $\mathbf{S}$, Eq.\,\eqref{eq:AdditiveSchwarz}. This effect can be resolved by the use of graph coloring algorithms \cite{Gebremedhin2005}.

\subsection{Loading of a die cast gearbox housing}\label{sec::gearbox}
Finally we consider an industrial example to portray the suitability of the preconditioner for the solution of large-scale embedded problems involving multi-level $hp$-refinement. The example focuses on the loading of an aluminum gearbox housing produced by die casting, a manufacturing process by which molten metal is forced under high pressure into a mold cavity producing a net-shape structure upon solidification of the molten metal. Die casting is widely used for the production of small- to medium-sized non-ferrous metal parts due to its high production rates, the quality and dimensional consistency of the produced parts and its design freedom, that allows the production of parts with high geometric complexity \cite{vinarcik2002high}. The extensive use of die cast parts is, however, limited by defects such as shrinkage cavities and entrapped gas bubbles that lead to increased porosity in the parts. Pores negatively affect the mechanical properties of die cast parts and are the preferred sites for fatigue-crack initiation.

Several studies have been conducted on the influence of pore size, position and orientation on the stress state, crack growth and fatigue behavior of aluminum die cast parts \citep{AMMAR20081024,MAYER2003245,Yi2003}. Finite element analysis of die cast parts can help quantify the stress concentration in the vicinity of pores through the use of image-based computations \cite{NICOLETTO2010547}. The limitations of boundary conforming methods in capturing the complex pore morphology can be overcome by the use of the finite cell method as investigated in \cite{monavari2011,Duczek2015,Wuerkner2018}. 

In this section, we revisit the die cast example considered in \cite{monavari2011}. The first major modification in contrast to the original example is the use of a preconditioned iterative solver for solving the linear system of equations instead of a direct solver. This modification allows, due to the smaller memory requirement, for computations with a higher resolution of the die cast part. Moreover, distributed memory parallelism is exploited to accelerate the performance of the preconditioned iterative solver. The second major modification is the use of multi-level $hp$-refinement to achieve a higher spatial resolution around the pores in the specimen leading to a better resolution of the complex stress state in these regions.    

\subsubsection*{Problem setup}\label{sec::gearboxsetup}
\begin{figure}[H]
 \begin{center}
  \begin{minipage}[c]{0.5\textwidth}
   \begin{center}
    \includegraphics[width=0.95\textwidth]{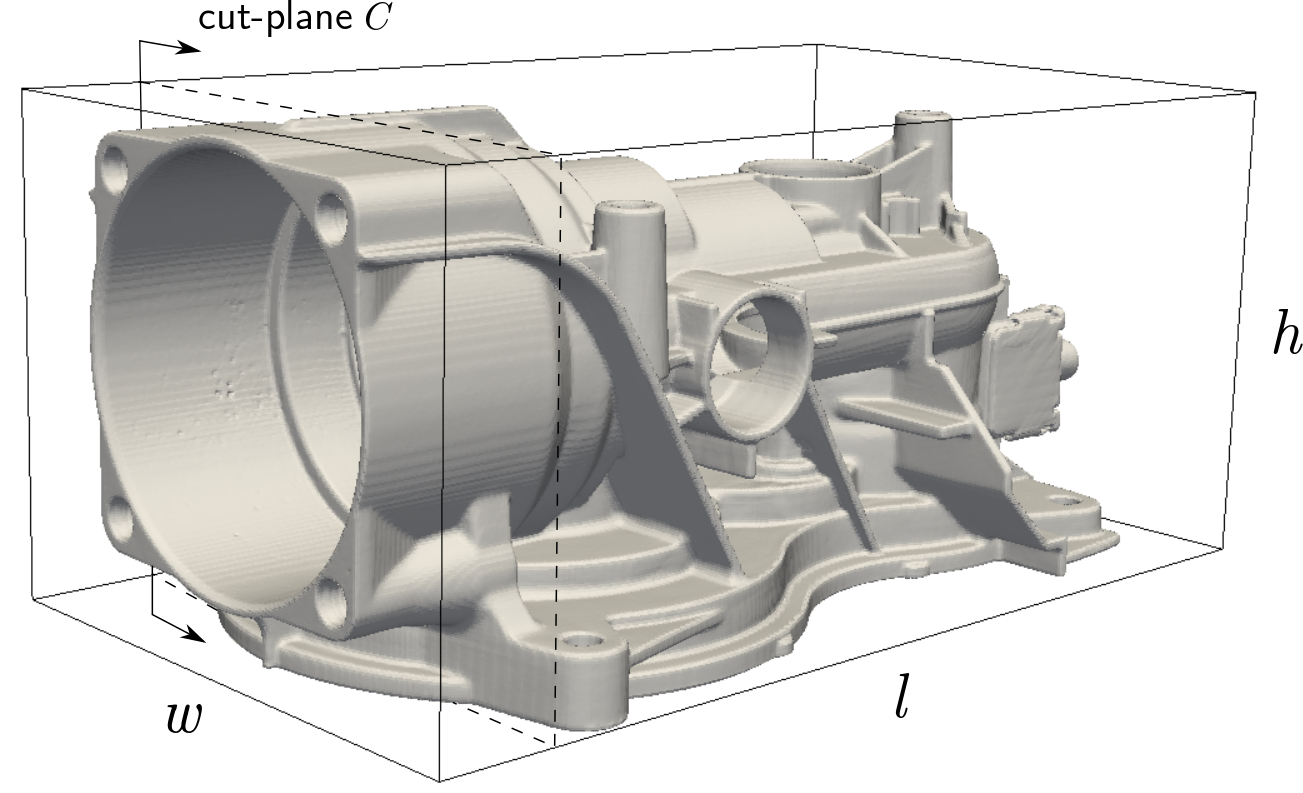}%
   \begin{table}[H]
   \centering
   \scriptsize
   \begin{tabular}{cc}
   \toprule
   \multirow{2}{*}{Material}  & aluminum alloy \\ & EN AC-46000 \\
   Young's modulus, $E$& 70 GPa \\
   Poisson's ratio, $\nu$ & 0.29 \\
   $l \times w \times h$ [mm] & $175 \times 101 \times 80 $   \\
   \bottomrule
   \end{tabular}
   \end{table}
   
   (a) Geometric and material properties of the gearbox housing. 

   \end{center}
  \end{minipage}
 \hspace{2mm}
 \begin{minipage}[c]{0.4\textwidth}
  \scriptsize
  \begin{center}
   \hspace*{1cm}\includegraphics[width=0.8\textwidth]{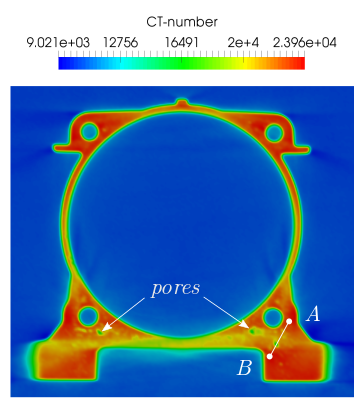}%
   \\ \hspace{0.5cm} (b) Cut-plane C of the housing CT scan 
   \\ \hspace{1.0cm} depicting pores in the specimen.
   \\
   \vspace{3mm}
   \includegraphics[width=0.8\textwidth]{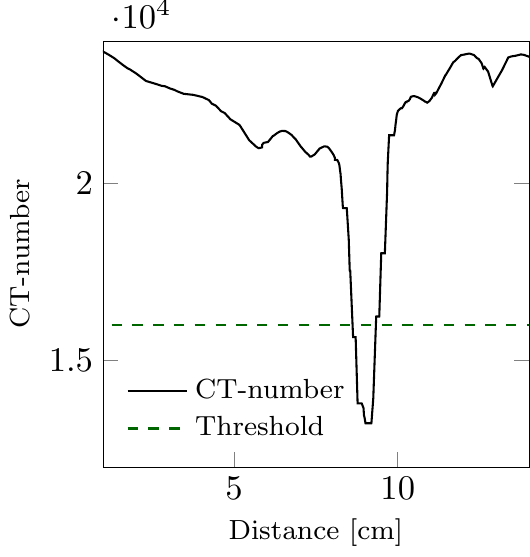}%
   \\ (c) Distribution of the CT-number \\ along the cut-line $A$-$B$. 
  \end{center}
  \end{minipage}
  \caption{Setup of the die cast gearbox housing example.}
  \label{fig::gearboxhousing}
 \end{center}
\end{figure}

\begin{figure}[H]
   \begin{center}
	   \subfloat[Discretization of the gearbox housing.]
     {
       \includegraphics[width=0.38\textwidth]{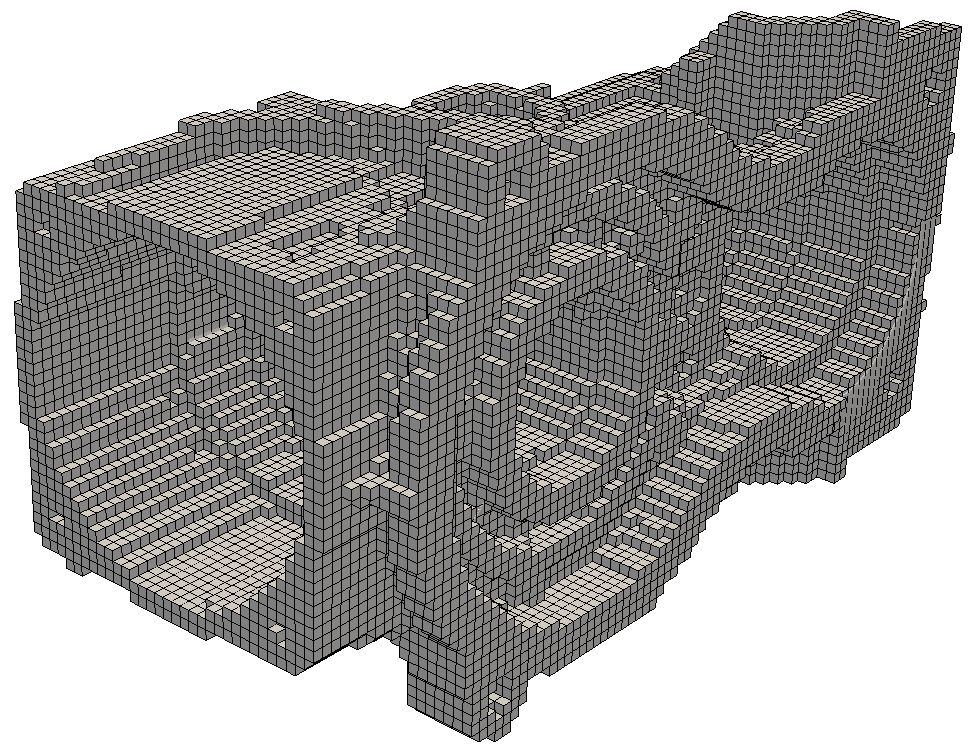}%
       \label{fig::gearboxhousing-mesh}
     }%
     \hspace{1cm}%
	   \subfloat[Application of boundary conditions.]
     {
       \includegraphics[width=0.43\textwidth]{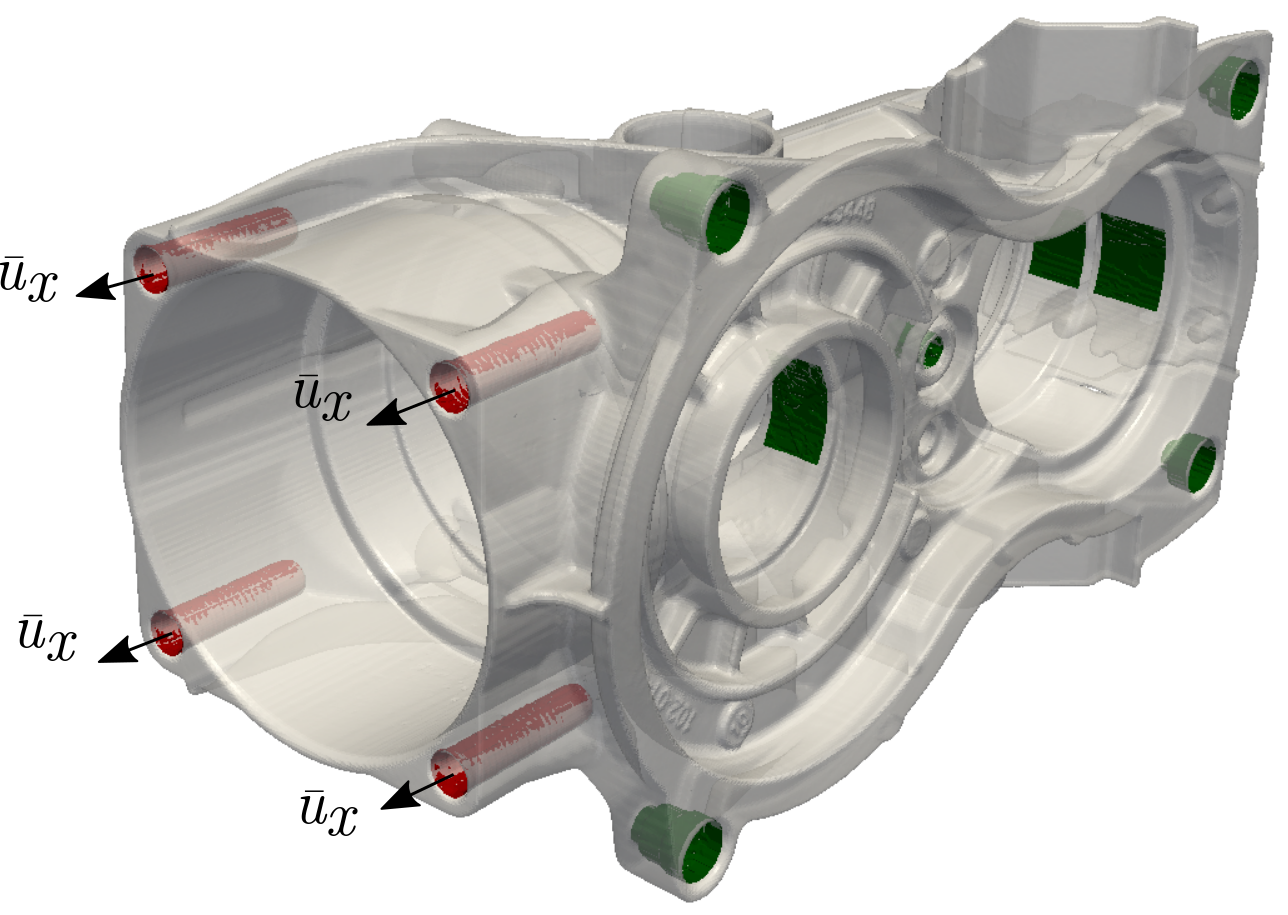}%
       \label{fig::gearboxhousing-boundaryconditions}
     }%
    \hfill%
    \end{center}
	\caption{Mesh and boundary conditions of the gearbox housing example. All displacements on the dark green surfaces are fixed while a displacement field $\bar{u}_{x}=0.1$cm is prescribed on the red surfaces.}
\end{figure}

The geometry of the gearbox housing is obtained from a CT scan with a voxel size of $198.7 \times 198.7 \times 198.7 \ \mu m^3$. A non-boundary conforming mesh is generated with $10^3$ voxels grouped to form a single finite cell. Elements completely inside the fictitious domain are eliminated yielding a total of 34\,230 active elements in the base mesh before refinement. The classification of the material properties is done by applying a threshold on the intensity values (CT-number) of the CT scan. The threshold is chosen such that the metal material and cavities or internal pores (which are part of the fictitious domain) can be distinguished, as illustrated in Figure \hyperref[fig::gearboxhousing]{\ref*{fig::gearboxhousing}b} and Figure \hyperref[fig::gearboxhousing]{\ref*{fig::gearboxhousing}c}\,. Metal material is assigned a Young's modulus $E_{\text{metal}}=70$ GPa and Poisson's ratio $\nu=0.3$ while fictitious material is assigned the values $E_{\text{fict}}=10^{-11}$ GPa and $\nu=0.3$. Figure \ref{fig::gearboxhousing-mesh} depicts the computational mesh of the gearbox housing. The discretization is refined towards the surfaces of the interior pores with a refinement depth $k \in \lbrace  0, 1, 2, 3 \rbrace$ and a polynomial order $p=4$ is chosen, resulting in the number of DOFs listed in Table \ref{tab:gearboxdofs}\,. 
The numerical integration in this example is is performed by Gauss quadrature over integration sub-cells instead of standard octrees, see e.g.\ \cite{Yang2012a,Elhaddad2017} for details.
Homogeneous Dirichlet boundary conditions are applied on the dark green surfaces while a displacement $\bar{u}_x=0.1$cm is applied on the red cylindrical surfaces as depicted in Figure \ref{fig::gearboxhousing-boundaryconditions}\,. All boundary conditions are enforced weakly using the penalty method with the penalty parameter $\beta=10^6$.

\begin{table}[H]
\centering
\begin{tabular}{ccccc}
\toprule
$k$ &  0  & 1 & 2 & 3 \\
\toprule
DOFs & 1\,662\,666 & 1\,670\,478  & 1\,683\,792  & 1\,716\,102 \\
\bottomrule
\end{tabular}
\caption{Number of degrees of freedom for different refinement depths $k$. }
\label{tab:gearboxdofs}
\end{table}

\subsubsection{Influence of pores on the stress distribution}

\begin{figure}[H]
  \centering
	\subfloat[Displacement magnitude.]
     {
       \includegraphics[width=0.45\textwidth]{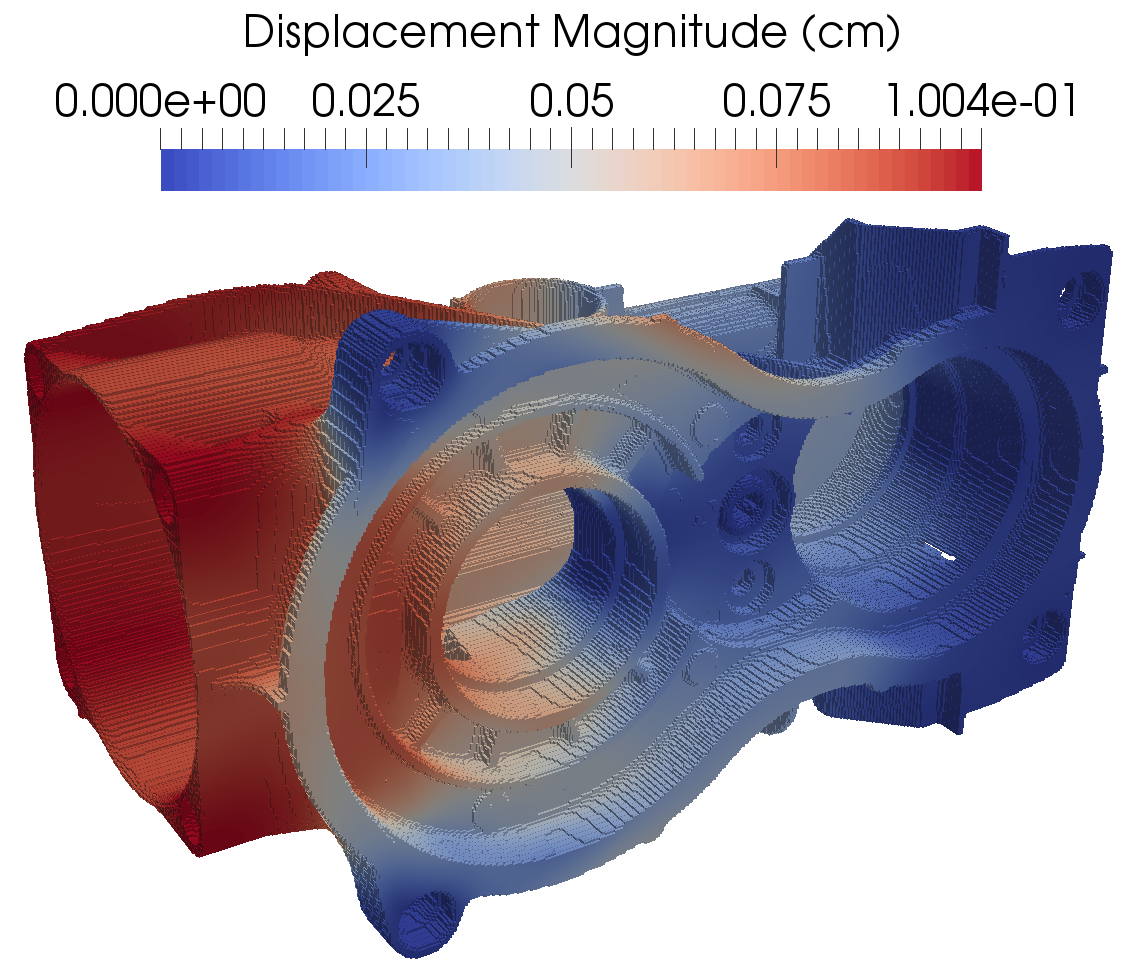}%
     }%
     \hfill%
	   \subfloat[Von Mises stress distribution.]
     {
       \includegraphics[width=0.45\textwidth]{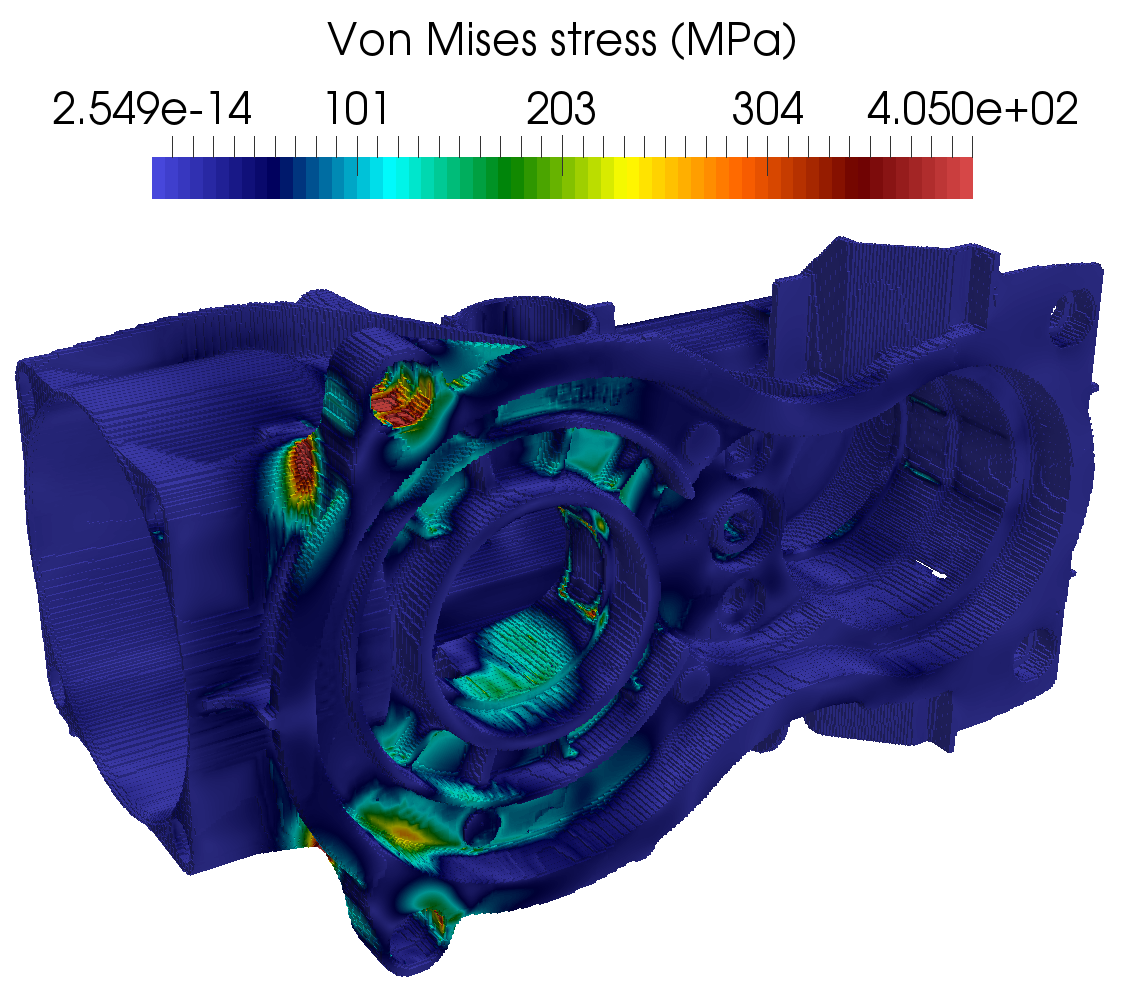}%
     }%
	\caption{Results of the displacement and von Mises stress fields of the gearbox housing.}
        \label{fig::globalresultsgearbox}
\end{figure}

\begin{figure}[H]
   \begin{center}
     {
       \includegraphics[width=0.45\textwidth]{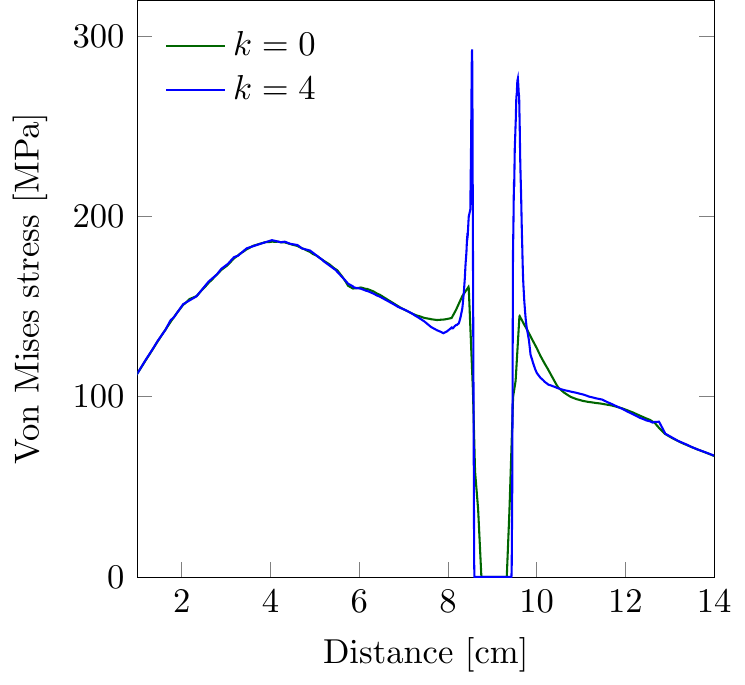}%
     }%
     \vspace{-1mm}
	   \caption{Comparison of the von Mises stress in the vicinity of the pores along the cut-line $A$-$B$ shown in Figure \hyperref[fig::gearboxhousing]{\ref*{fig::gearboxhousing}b}\,. }
        \label{fig::lineresultsgearbox}
    \end{center}
\end{figure}

\begin{figure}[H]
   \begin{center}
    \subfloat[Distribution of the von Mises stress in cut-plane $C$ for $k=0$.] 
    {       
            \label{fig::localresultsgearboxA}
            \includegraphics[]{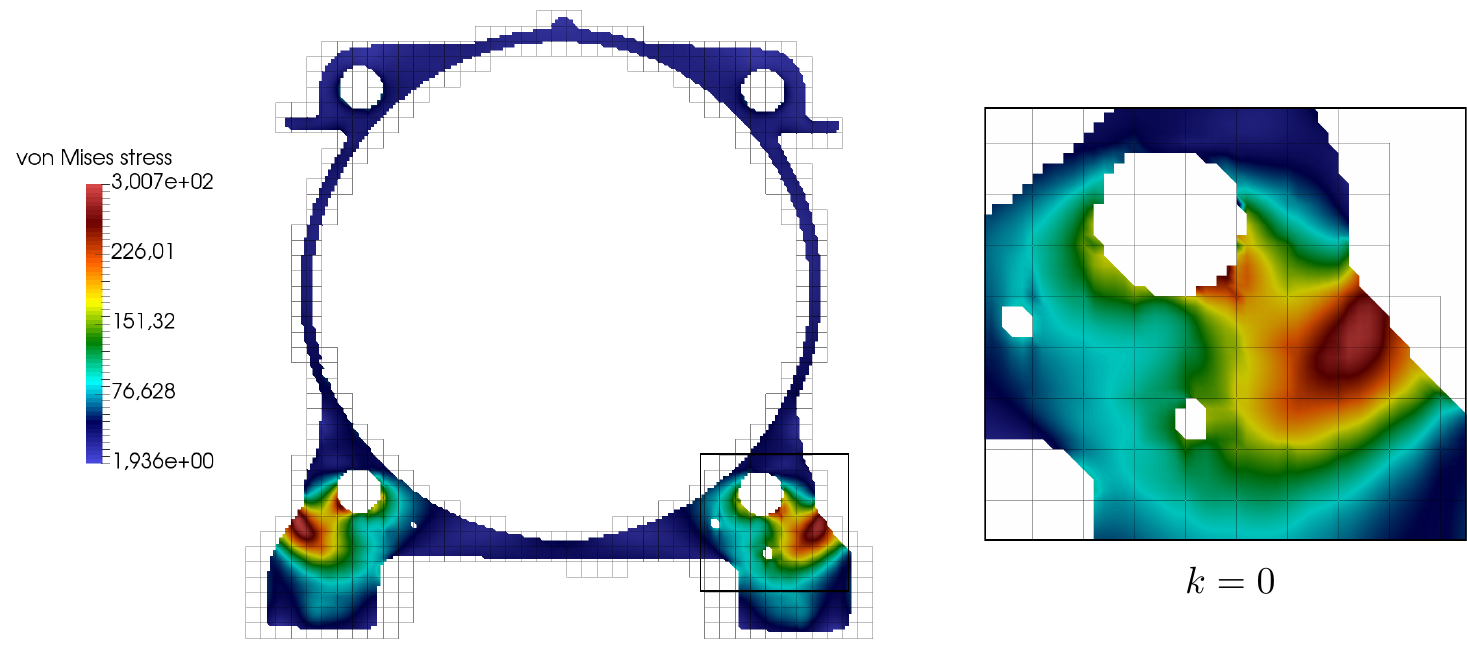}
    }
    \hfill
   \subfloat[Distribution of the von Mises stress in cut-plane $C$ for $k=4$.]
   {
         \label{fig::localresultsgearboxB}
         \includegraphics[]{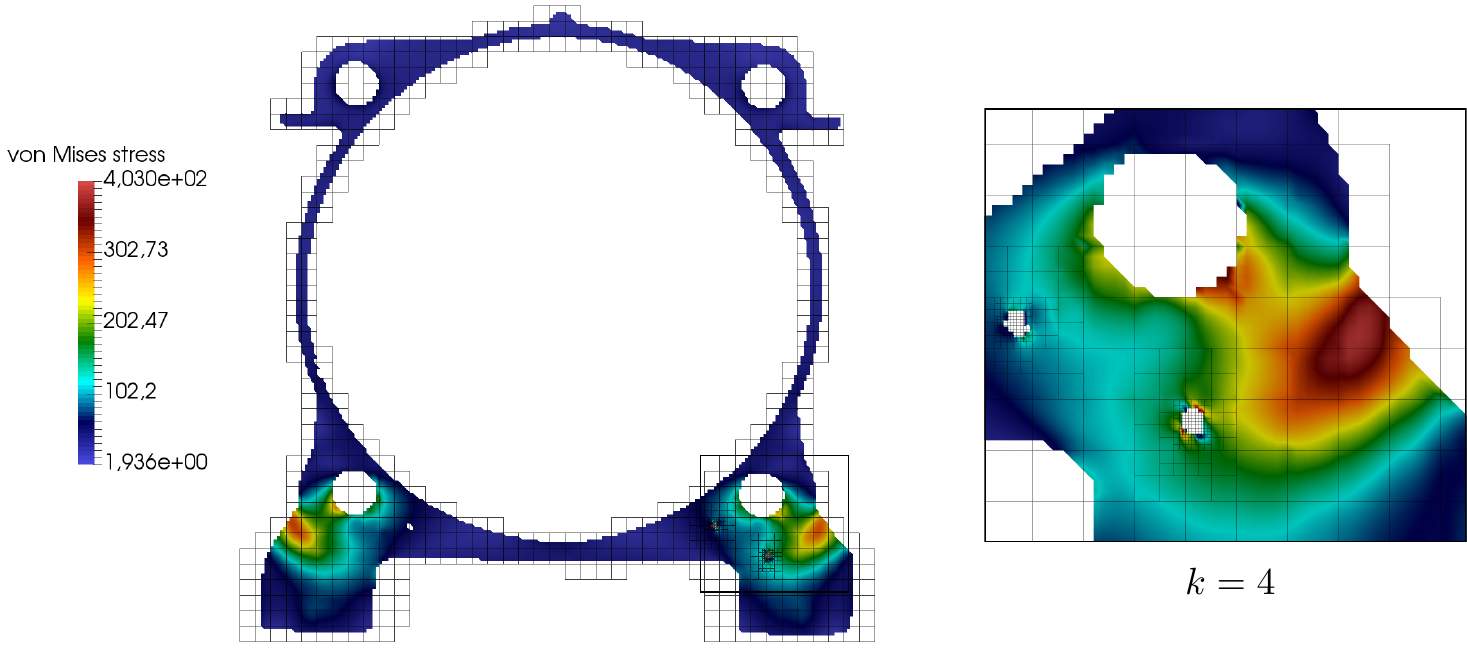}
   }
\caption{Influence of the refinement depth on the stress distribution around pores in cut-plane $C$.}
\label{fig::localresultsgearbox}
    \end{center}
\end{figure}

Before studying the convergence behavior of the iterative solver for different refinement depths $k$, we highlight the necessity of local refinements in resolving stress states in the vicinity of small geometric features with an economical number of DOFs. Different values of $k$ yield comparable results of the overall displacement and stress fields in the gearbox housing with an exemplary result given in Figure \ref{fig::globalresultsgearbox}\,. Differences arise, however, in the stress fields in the vicinity of the pores as shown in Figures \ref{fig::lineresultsgearbox} - \ref{fig::localresultsgearbox}\,. Without local mesh refinement, the chosen mesh is unable to accurately capture the expected stress concentration around the pore boundary, as shown in Figure \ref{fig::localresultsgearboxA}\,. The use of a globally $h$-refined mesh would help capture stress concentration around the pores, but is not a viable option since the number of DOFs would increase drastically. Moreover, performing $h$-refinement would not significantly affect the stress distribution away from the pores, since this state is already adequately represented by the coarse mesh as shown in Figure \ref{fig::lineresultsgearbox}\,. The stated arguments motivate the use of $hp$-refinement to capture local solution stress distributions in a computationally efficient manner. The results of the Von Mises stress available in Figure \ref{fig::localresultsgearboxB} and Figure \ref{fig::lineresultsgearbox} show that the location of the peak stress shifts to the boundary of the pores upon application of multi-level $hp$-refinement, which is physically realistic behavior. 

\subsubsection{Convergence behavior}
\begin{figure}[H]
   \begin{center}
	   \subfloat[Full preconditioner: residual convergence. ]
     {
       \includegraphics[width=0.48\textwidth]{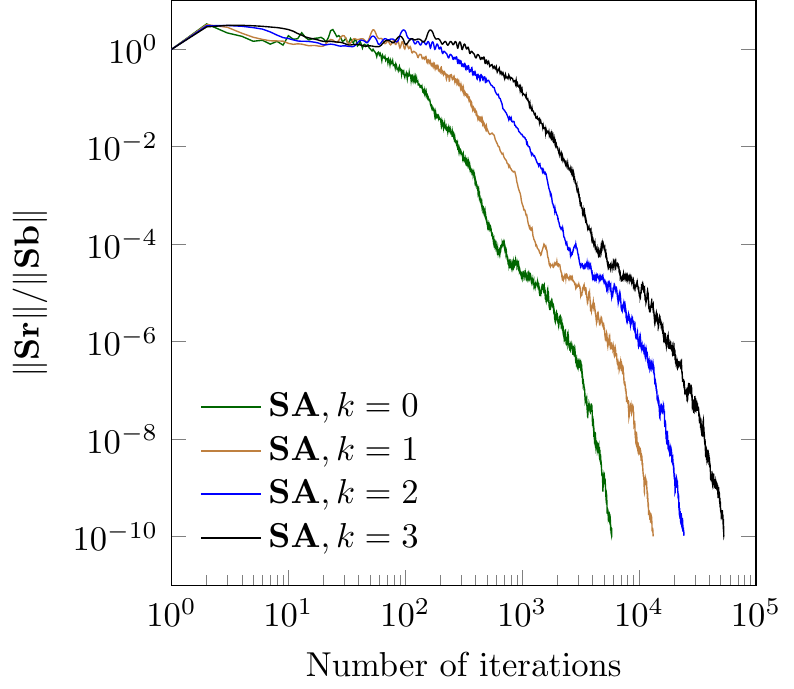}%
     }%
     \hfill%
	   \subfloat[Truncated preconditioner: residual convergence. ]
     {
       \includegraphics[width=0.48\textwidth]{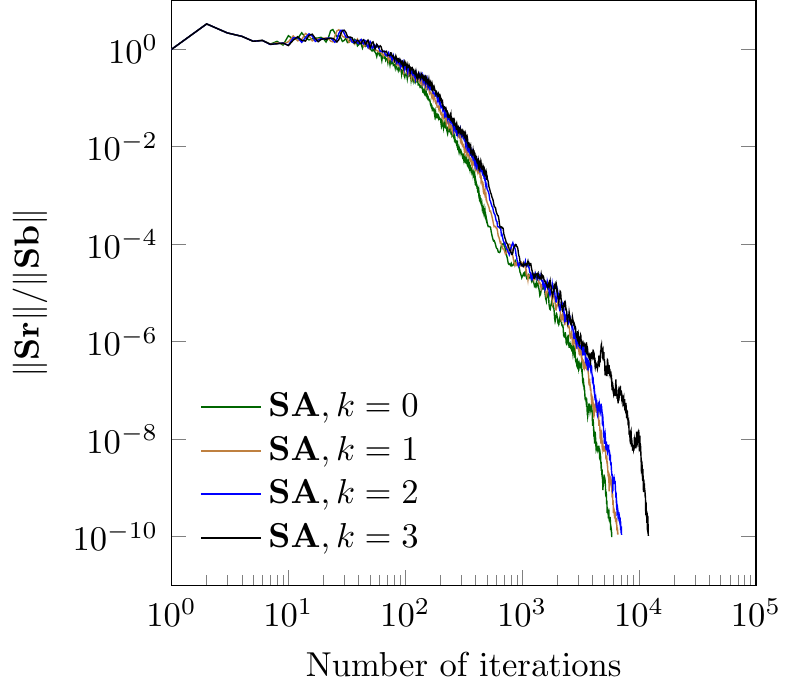}%
     }%
     \hfill%
	   \subfloat[Full preconditioner: energy convergence. ]
     {
       \includegraphics[width=0.48\textwidth]{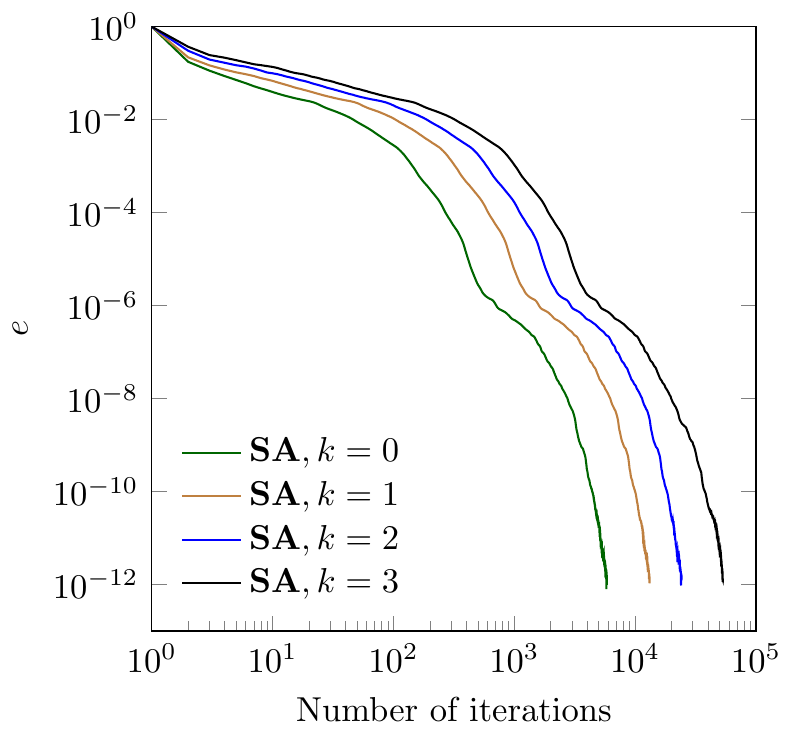}%
     }%
     \hfill%
	   \subfloat[Truncated preconditioner: energy convergence. ]
     {
       \includegraphics[width=0.48\textwidth]{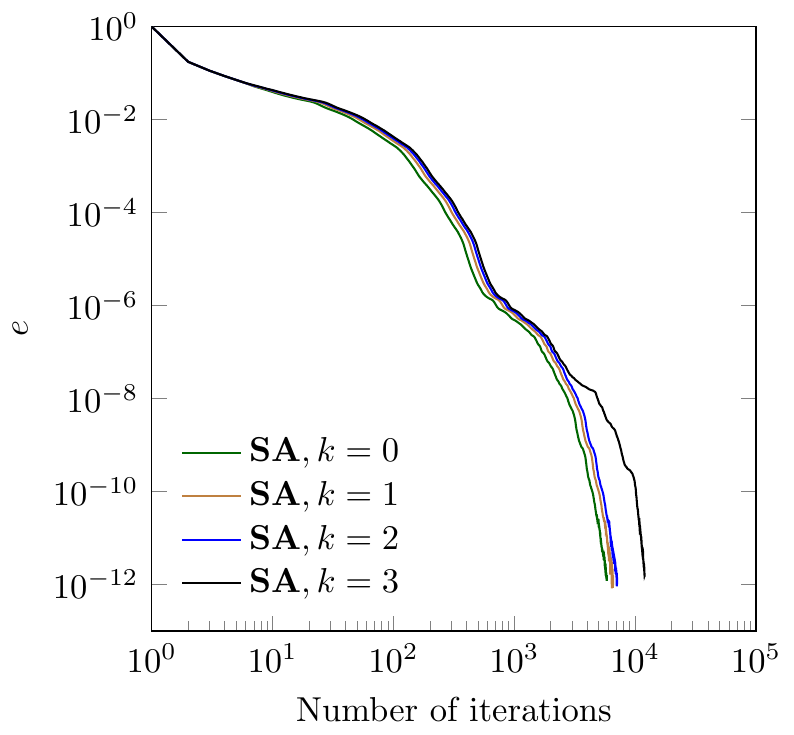}%
     }%
	\caption{ Comparison of the convergence behavior of the full and truncated preconditioners for multi-level $hp$-refinement.}
        \label{fig::convergencegearbox}
    \end{center}
\end{figure}

An analogous study to that in Section \ref{sec::ballinbox} is carried out to analyze the influence of multi-level $hp$-refinement on the convergence of the PCG solver. We refine the computational mesh towards the internal pores with a refinement depth $k \in \lbrace 0, 1, 2, 3 \rbrace $.  The development of the relative preconditioned residual and energy error is  monitored in computations involving the full and truncated preconditioners introduced in Section \ref{sec::preconditioningmlhprefinement}\,. Figure \ref{fig::convergencegearbox} portrays the convergence behavior of the full and truncated preconditioner for different levels of multi-level $hp$-refinement of the gearbox housing mesh. The results show similar convergence behavior as in the numerical example considered in Section \ref{sec::ballinbox}\,. When the full preconditioner is used, the number of iterations drastically increases with the refinement depth. The truncated preconditioner performs better than the full preconditioner since it only takes into account basis functions that can potentially become almost linearly dependent as described in Section \ref{sec::preconditioningmlhprefinement}\,. Table~\ref{tab:gearboxgroups} demonstrates that $n_{\textrm{max}}$, the maximum number of blocks a basis function belongs to, is bounded and equals $2^d=8$ for the truncated preconditioner.

\begin{table}[H]
\begin{center}
\begin{tabular}{ccccc}
\toprule
$k$ &  0  & 1 & 2 & 3 \\
\toprule
Full blocks & 8 & 64 & 274  & 810 \\
Truncated blocks & 8 & 8  & 8  & 8 \\
\bottomrule
\end{tabular}
\caption{Maximum overlap $n_{\textrm{max}}$ of the full and truncated Additive-Schwarz blocks for different refinement depths $k$ in the gearbox housing example. }
\label{tab:gearboxgroups}
\end{center}
\end{table}

\subsubsection{Analysis of the parallel scalability}\label{sec::mpiscalability}

Section \ref{sec::mpiparallelism} describes the use of the preconditioner in a hybrid (distributed and shared memory) parallel setting. The study at hand aims to show the performance of a parallel PCG solver in large-scale finite cell computations. We use the parallel PCG solvers available in \texttt{Trilinos} \cite{Trilinos} in conjunction with the presented preconditioner. Two different discretizations of the gearbox housing are considered. The first discretization is the mesh described in Section \ref{sec::gearboxsetup} consisting of 34\,230 elements and approximately 1.6 million degrees of freedom. The second mesh is a finer discretization of the gearbox that consists of 83\,757 elements formed by grouping $7^3$ voxels to form a single finite cell. A polynomial order of $p=4$ is chosen here as well, resulting in approximately 3.98 million DOFs. Boundary conditions are applied for both discretizations as illustrated in Figure \ref{fig::gearboxhousing-boundaryconditions} with a penalty value $\beta=10^{6}$ . The numerical computations are performed on the CoolMAC cluster at the Technical University of Munich equipped with four AMD Bulldozer Opteron 6274 CPUs and 256 GB memory per node. These hybrid computations are carried out by varying the number of nodes from 1 to 16 for the coarser discretization and from 2 to 16 for the fine discretization. Each node is made up of 64 processing units (cores) i.e. one MPI-process $\times$ 64 OpenMP threads.  

\begin{figure}[H]
   \begin{center}
       \includegraphics[width=0.48\textwidth]{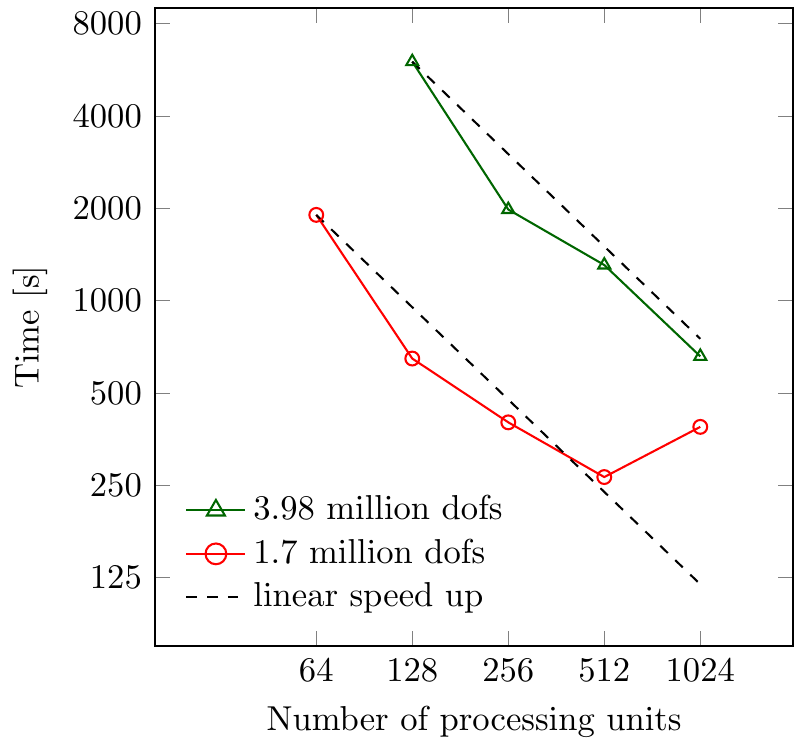}%
	\caption{Parallel scalability study of the PCG solver for two different problem sizes.}
        \label{fig::hybridscalingAS}
    \end{center}
\end{figure}

Figure \ref{fig::hybridscalingAS} shows the execution time of the PCG solver plotted against the number of processing units used in the hybrid computations for both discretizations of the gearbox housing. The preconditioned solver shows excellent parallel scalability for both discretizations up to 256 cores where superlinear speed up is achieved due to cache effects when moving to more cores. This trend can be maintained as long as the computational work in every PCG iteration --- the matrix vector multiplications --- is significantly larger than the communication overhead. This is shown by the loss of efficiency in the computations involving 1.7 million degrees of freedom and more than 256 cores. The fine discretization provides adequate computational work and maintains a parallel efficiency of 115\% up to 1024 processing units.

%% file: conclusion/conclusion.tex
\section{Conclusion}

The presented work describes a preconditioning technique that enables iterative solutions of large finite cell method (FCM) computations with multi-level $hp$-refined discretizations. An Additive Schwarz preconditioning technique based on overlapping element-wise blocks is developed and tailored to the hierarchical structure of the $hp$-discretization scheme. Furthermore, a robust preconditioner is obtained through the use of a pseudo inverse and the efficiency of the preconditioner improved by utilizing parallelism. The results show that the preconditioner is robust to local refinements and to how elements are cut. It is verified that the computational cost of setting up and applying the preconditioner scales well with the number of processors in a parallel setting. Therefore the preconditioning technique increases tolerable sizes of FCM systems, which extends the scope of problems that can be solved with FCM to real-life applications.

The current work is restricted to problems which yield symmetric positive definite matrices. Notwithstanding that this covers a large range of engineering problems, this does not include multivariate problems with indefinite matrices or problems involving convection with nonsymmetric matrices. Applying or modifying the preconditioner for such applications, as done for simple meshes in \cite{CbAS}, is therefore an important topic of future research. Furthermore, the current setup is tailored for the multi-level $hp$-refinement scheme introduced in \cite{Zander2015}. A relevant topic for further research is the development and testing of similar techniques for different refinement schemes.

Finally the preconditioner does not resolve the conditioning effects that can be expected from increasing the number of elements, e.g.\ \cite{johnson2012numerical}. Be it that the preconditioner is robust to cut elements and thereby extends the numbers of DOFs that FCM systems can currently be solved with, we do not expect it to be efficient for exascale systems. In order to further extend the range of solvable systems, it is an interesting topic for future research to apply the preconditioner together with a multi-grid or domain-decomposition approach.

%% file: template/acknowledgements.tex
The authors gratefully acknowledge the financial support of the NWO under the Graduate program Fluid \& Solid Mechanics and the European Research Council under Grant ERC-2014-
StG 637164.